# A cardiac electromechanics model coupled with a lumped parameters model for closed-loop blood circulation. Part II: numerical approximation

Francesco Regazzoni[1,*], Matteo Salvador[1,*], Pasquale Claudio Africa[1], Marco Fedele[1], Luca Dedè[1], Alfio Quarteroni[1,2]

[1] MOX-Dipartimento di Matematica, Politecnico di Milano, Milan, Italy
[2] Professor Emeritus, École polytechnique fédérale de Lausanne, Lausanne, Switzerland

∗ These authors equally contributed to this work

**Abstract**

In the framework of accurate and efficient segregated schemes for 3D cardiac electromechanics and 0D cardiovascular models, we propose here a novel numerical approach to address the coupled 3D-0D problem introduced in Part I of this two-part series of papers. We combine implicit-explicit schemes to solve the different cardiac models in a multiphysics setting. We properly separate and manage the different time and space scales related to cardiac electromechanics and blood circulation. We employ a flexible and scalable intergrid transfer operator that enables to interpolate Finite Element functions among different meshes and, possibly, among different Finite Element spaces. We propose a numerical method to couple the 3D electromechanical model and the 0D circulation model in a numerically stable manner within a fully segregated fashion. No adaptations are required through the different phases of the heartbeat. We also propose a robust algorithm to reconstruct the stress-free reference configuration. Due to the computational cost associated with the numerical solution of this inverse problem, the reference configuration recovery algorithm comes along with a novel projection technique to precisely recover the unloaded geometry from a coarser representation of the computational domain. We show the convergence property of our numerical schemes by performing an accuracy study through grid refinement. To prove the biophysical accuracy of our computational model, we also address different scenarios of clinical interest in our numerical simulations by varying preload, afterload and contractility. Indeed, we simulate physiologically relevant behaviors and we reproduce meaningful results in the context of cardiac function.

## 1 Introduction

This paper proposes novel numerical approaches for the simulation of cardiac electromechanics [12, 16, 20, 35, 44, 54, 58]. Following [49], we propose the numerical approximation and simulation of the 3D-0D model provided in (1). The computational model for cardiac electromechanics presented in this second part relies on segregated schemes [17, 19, 21]. These schemes, also called partitioned, often suffer from instabilities issues [19, 39, 43, 50], despite being computationally more attractive than monolithic schemes [21, 22, 59]. Nonetheless, our numerical approach enable to couple the core models describing the different physics in a numerically stable manner, yet allowing to adopt different resolutions in space and time for each one, to reflect their characteristic scales [9, 20, 35, 44, 52, 56]. With this aim, we implement an efficient and flexible intergrid transfer operator [1, 23, 52] that allows to couple models defined on different computational meshes and possibly different Finite Element spaces. The single core models are approximated by suitable implicit-explicit (IMEX) schemes [42] in time, aimed at lowering the computational effort associated with each time step without introducing severe CFL restrictions on the time step size otherwise associated to fully explicit schemes [45]. The computational framework presented in this paper realizes a very favorable trade-off between the biophysical detail of the underlying mathematical models and computational efficiency of the corresponding solver. Mathematical rigor, model accuracy and computational efficiency are the landmarks of the proposed electromechanical solver.

We also propose a novel numerical approach to couple the 3D and the 0D models [8, 28] in the segregated scheme. We reinterpret the cavity pressure of the left ventricle as a Lagrange multiplier to enforce a volumetric constraint to couple together the 0D circulation model and the 3D electromechanics model. We end up with a saddle-point structure for the mechanical problem, which is numerically solved using the Schur complement reduction [7]. Our scheme is numerically stable and can be applied to all the phases of the heart cycle (filling phase, ejection phase and isovolumetric contraction and relaxation phases) without the need of switching among different 0D models and without considering any change in the parameters of the equations [34, 52].

Our model is supported by a strategy to reconstruct the stress-free (i.e. unloaded) configuration of a ventricle, starting from the loaded geometry acquired by cardiac imaging, by solving a suitable inverse problem [29]. This is essential for physically consistent patient-specific simulations since the elastodynamics equations must refer to such stress-free reference configuration. The numerical resolution of this inverse problem can be computationally demanding for highly refined meshes. For this reason, our algorithm for the recovery of the reference configuration is supplied by a projection technique that accurately recover the stress-free configuration from a coarser and independent representation of the computational domain.

In terms of numerical simulations, we perform a grid refinement study, using also different Finite Element spaces among the different core models, to show the effectiveness of our numerical approach. Moreover, we simulate the response of our electromechanical model to some physiologically relevant scenarios, namely varying preload, afterload and contractility.

### 1.1 Paper outline

In Sec. 2 we briefly recall the mathematical model of cardiac electromechanics considered in this paper. In Sec. 3 we present the numerical discretization of the different core models and, more importantly, the strategy followed for their coupling. In Sec. 4 we illustrate our algorithm to recover the reference configuration from a stressed geometry. In Sec. 5 we present the numerical results: we perform a grid refinement study and we analyze different numerical tests by varying preload, afterload and contractility, thus showing the suitability of our model to investigate different physiologically relevant scenarios. Finally, in Sec. 6, we draw our conclusions and outlook on future developments.

## 2 Mathematical models

We recall the models proposed in the Part I of this paper [49]:

$$
(\mathscr{E})\begin{cases}\chi_\mathrm{m}\left[C_\mathrm{m}\dfrac{\partial u}{\partial t}+\mathcal{I}_\mathrm{ion}(u,\boldsymbol{w},\boldsymbol{z})\right]-\nabla\cdot\left(J\mathbf{F}^{-1}\boldsymbol{D}_\mathrm{M}\mathbf{F}^{-T}\nabla u\right)\\ \qquad\qquad\qquad\qquad\qquad\qquad\qquad\qquad\qquad=\mathcal{I}_\mathrm{app}(t),\\ \left(J\mathbf{F}^{-1}\boldsymbol{D}_\mathrm{M}\mathbf{F}^{-T}\nabla u\right)\cdot\mathbf{N}=0,\end{cases} \tag{1a}
$$

in $\Omega_0\times(0,T)$, with $u=u_0$ in $\Omega_0$ at time $t=0$.

$$
(\mathscr{I})\begin{cases}\dfrac{\partial\boldsymbol{w}}{\partial t}-\boldsymbol{H}(u,\boldsymbol{w})=\mathbf{0},\\ \dfrac{\partial\boldsymbol{z}}{\partial t}-\boldsymbol{G}(u,\boldsymbol{w},\boldsymbol{z})=\mathbf{0},\end{cases} \tag{1b}
$$

in $\Omega_0\times(0,T)$, with $\boldsymbol{w}=\boldsymbol{w}_0$ and $\boldsymbol{z}=\boldsymbol{z}_0$ in $\Omega_0$ at time $t=0$.

$$
(\mathscr{A})\quad\frac{\partial\mathbf{s}}{\partial t}=\boldsymbol{K}(\mathbf{s},[\mathrm{Ca}^{2+}]_\mathrm{i},SL), \tag{1c}
$$

in $\Omega_0\times(0,T)$, with $\mathbf{s}(0)=\mathbf{s}_0$ in $\Omega_0$ at time $t=0$.

$$
(\mathscr{M})\begin{cases}\rho_\mathrm{s}\dfrac{\partial^2\mathbf{d}}{\partial t^2}-\nabla\cdot\mathbf{P}(\mathbf{d},T_\mathrm{a}(\mathbf{s}))=\mathbf{0},\\ \mathbf{P}(\mathbf{d},T_\mathrm{a}(\mathbf{s}))\mathbf{N}+\mathbf{K}^\mathrm{epi}\mathbf{d}+\mathbf{C}^\mathrm{epi}\dfrac{\partial\mathbf{d}}{\partial t}=\mathbf{0} & \text{on }\Gamma_0^\mathrm{epi}\times(0,T),\\ \mathbf{P}(\mathbf{d},T_\mathrm{a}(\mathbf{s}))\mathbf{N}=p_\mathrm{LV}(t)\,|J\mathbf{F}^{-T}\mathbf{N}|\,\mathbf{v}^\mathrm{base}(t) & \text{on }\Gamma_0^\mathrm{base}\times(0,T),\\ \mathbf{P}(\mathbf{d},T_\mathrm{a}(\mathbf{s}))\mathbf{N}=-p_\mathrm{LV}(t)\,J\mathbf{F}^{-T}\mathbf{N} & \text{on }\Gamma_0^\mathrm{endo}\times(0,T),\end{cases} \tag{1d}
$$

in $\Omega_0\times(0,T)$, with $\mathbf{d}=\mathbf{d}_0$ and $\dfrac{\partial\mathbf{d}}{\partial t}=\dot{\mathbf{d}}_0$ in $\Omega_0$ at time $t=0$.

$$
(\mathscr{C})\quad\frac{d\boldsymbol{c}_1(t)}{dt}=\widetilde{\boldsymbol{D}}(t,\boldsymbol{c}_1(t),p_\mathrm{LV}(t)), \tag{1e}
$$

for $t\in(0,T)$, with $\boldsymbol{c}_1(0)=\boldsymbol{c}_{1,0}$.

$$
(\mathscr{V})\quad V_\mathrm{LV}^\mathrm{0D}(\boldsymbol{c}_1(t))=V_\mathrm{LV}^\mathrm{3D}(\mathbf{d}(t)), \tag{1f}
$$

for $t\in(0,T)$.

In these models, $u\colon\Omega_0\times(0,T)\to\mathbb{R}$ denotes the transmembrane potential, $\boldsymbol{w}\colon\Omega_0\times(0,T)\to\mathbb{R}^{n_{\boldsymbol{w}}}$ and $\boldsymbol{z}\colon\Omega_0\times(0,T)\to\mathbb{R}^{n_{\boldsymbol{z}}}$ the ionic variables (respectively describing the opening of ionic channels and the ionic concentrations). The vector-valued variable $\mathbf{s}\colon\Omega_0\times(0,T)\to\mathbb{R}^{n_\mathbf{s}}$ collects the states variables of the force generation model. The variable $\mathbf{d}\colon\Omega_0\times(0,T)\to\mathbb{R}^3$ denotes the displacement of the tissue. Finally, $\boldsymbol{c}_1\colon(0,T)\to\mathbb{R}^{n_\mathbf{c}}$ is the state vector of the circulation models (collecting the pressures, volumes and fluxes in the different compartments of the vascular network) and $p_\mathrm{LV}\colon(0,T)\to\mathbb{R}$ denotes the left ventricle pressure.

In Eq. (1), cardiac electrophysiology is described by the monodomain equation $(\mathscr{E})$ [14], coupled with the ten Tusscher-Panfilov ionic model $(\mathscr{I})$ (henceforth denoted by TTP06) [55]. The generation of active force $(\mathscr{A})$ is described by means of an Artificial Neural Network (ANN) based model that surrogates the high-fidelity model proposed in [46] (RDQ18 model). This reduced model, built by exploiting the Machine Learning algorithm proposed in [47, 48], drastically reduces the computational cost associated with the resolution of the RDQ18 model, still reproducing its results with a very good accuracy [48]. In ($\mathscr{M}$), we model the passive elastic properties of the cardiac tissue through the Guccione constitutive model [26, 27]. Finally, blood circulation is described by means of the lumped parameters model introduced in Part I of this paper [49].

# 3 Numerical approximation

We sketch in Fig. 1 the numerical approximation of problem (1). We use a segregated numerical scheme that allows to separate and properly manage the space and time scales related to cardiac electromechanics, by using different mesh sizes and different time steps according to the characteristic scale of each core model. Our scheme is thus segregated (i.e. the different core models are solved sequentially) and staggered (i.e. possibly different time step sizes are used for the different core models) [21, 52].

We perform the numerical approximation of the core models using the Finite Element Method (FEM) in space and Finite Difference schemes in time [42]. We employ an explicit $4^{th}$ order Runge-Kutta method for the set of ODEs ($\mathscr{C}$). This high-order time discretization is motivated by the stiffness of the circulation model [42].

We consider a fine representation of the computational domain for both ($\mathscr{I}$) and ($\mathscr{E}$) models, whereas a coarser one is employed for both ($\mathscr{A}$) and ($\mathscr{M}$). This is motivated by the requirement of a higher resolution for ($\mathscr{I}$) and ($\mathscr{E}$), due to the sharp wavefronts characterizing electrophysiological solutions, whereas both ($\mathscr{A}$) and ($\mathscr{M}$) feature larger spatial scales [13, 15, 52]. Moreover, the nonlinearities of ($\mathscr{M}$) warrants for the use of a coarser space discretization to make the numerical solution less computationally intensive.

For what concerns time discretization, we use a staggered strategy based on the Godunov splitting scheme [24]. This approach introduces a first order splitting error [17, 21], which is compatible with the first order error associated with the backward Euler method used for ($\mathscr{I}$), ($\mathscr{E}$), ($\mathscr{A}$) and ($\mathscr{M}$)–($\mathscr{V}$). Moreover, the Godunov splitting scheme does not undermine the stability of the numerical scheme [17]. In our pipeline, we first update the variables of ($\mathscr{I}$) and ($\mathscr{E}$), then the variables of ($\mathscr{A}$) and finally, after updating the unknowns of ($\mathscr{C}$), we update the ones of ($\mathscr{M}$)–($\mathscr{V}$).

For electrophysiology, we employ a semi-implicit scheme, that we denote by ($\mathscr{I}_{\mathrm{SI}}$)–($\mathscr{E}_{\mathrm{SI}}$), where SI stands for semi-implicit. Mechanical activation ($\mathscr{A}_{\mathrm{E}}$), is solved with an explicit method in time. Mechanics ($\mathscr{M}_{\mathrm{I}}$) is instead implicitly discretized in time, due to the fact that the highly nonlinear (exponential) terms of the strain energy function $\mathcal{W}$ would need a restriction on the time step in both the semi-implicit and explicit contexts. Finally, we employ an explicit $4^{th}$ order Runge-Kutta scheme for the circulation model, indicated as ($\mathscr{C}_{\mathrm{E}}$). We use two different time steps for ($\mathscr{I}_{\mathrm{SI}}$)–($\mathscr{E}_{\mathrm{SI}}$)–($\mathscr{A}_{\mathrm{E}}$)–($\mathscr{C}_{\mathrm{E}}$) and for ($\mathscr{M}_{\mathrm{I}}$)–($\mathscr{V}_{\mathrm{I}}$).

We propose a novel strategy to perform the 3D-0D coupling, given by ($\mathscr{V}$), between ($\mathscr{M}$) and ($\mathscr{C}$), at the numerical level, within a segregated setting. Specifically, we solve the 3D mechanical problem under a volumetric constraint coming from the 0D circulation model. The cavity pressure acts as Lagrange multiplier associated with this constraint. Thus, we obtain a saddle-point problem that we address, at the algebraic level, by means of a Schur complement reduction.

## 3.1 Space discretization and intergrid interpolation

We implemented an intergrid transfer operator between nested grids that can be generalized to the case of locally-refined non-conforming nested grids. More details about this strategy can be found in [1]. We built an efficient and scalable interpolation data structure on top of the intergrid transfer operator that enables to evaluate the feedback coming from ($\mathscr{E}$) at the quadrature nodes defined on the mesh used for ($\mathscr{M}$) (and vice-versa), regardless of the polynomial degree used to discretize the two different models. This enables for a high level of numerical flexibility: different Finite Element degrees, as well as different levels of mesh refinement, can be effortlessly selected to tune the computational efficiency on the desired accuracy.

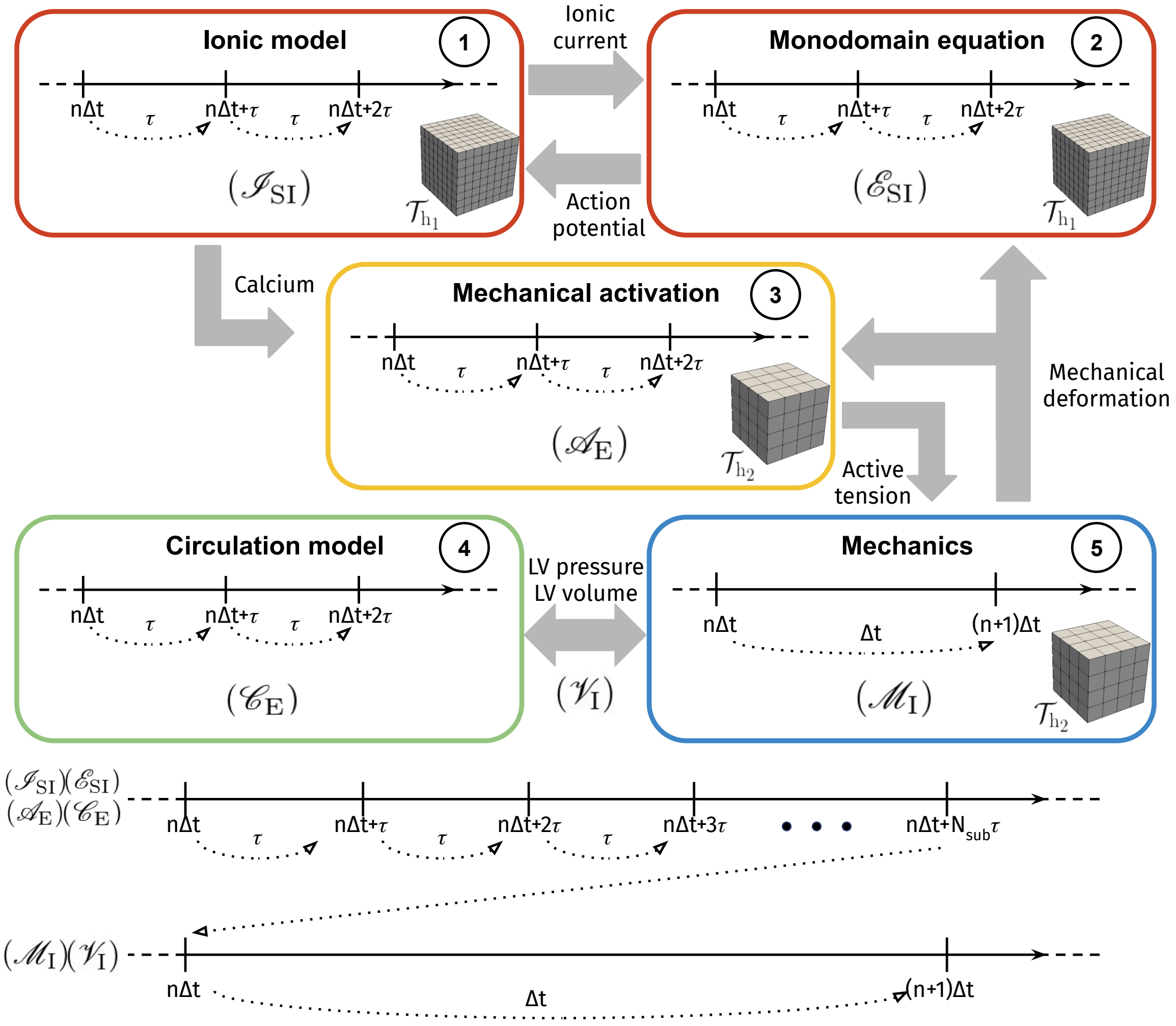

Figure 1: Sketch of the segregated-intergrid-staggered numerical scheme.

We consider two nested hexahedral meshes $\mathcal{T}_{\mathrm{h}_1}$ and $\mathcal{T}_{\mathrm{h}_2}$ of the computational domain $\Omega_0$, where $\mathcal{T}_{\mathrm{h}_1}$ has been generated by uniformly refining $\mathcal{T}_{\mathrm{h}_2}$ according to an octree structure [1, 11], i.e. by recursively splitting each parent element of $\mathcal{T}_{\mathrm{h}_2}$ into eight sub-elements for a prescribed number of times, i.e. until the desired geometrical detail is reached. Here $\mathrm{h}_1$ and $\mathrm{h}_2$ (with $\mathrm{h}_1 < \mathrm{h}_2$) represent the mesh sizes, which we compute as the mean of the maximum diameter of each element.

We denote by $N_{\mathrm{u}}$, $N_{\boldsymbol{w},\boldsymbol{z}}$, $N_{\mathbf{s}}$ and $N_{\mathbf{d}}$ the number of degrees of freedom (DOFs, that is the number of variables) for the dimensionless transmembrane potential, gating and concentration variables, mechanical activation variables and displacement, respectively. We denote the set of tensor-products of polynomials with degree smaller than or equal to $r$ over a mesh element $K$ by $\mathbb{Q}_{\mathrm{r}}(K)$, and we introduce the finite dimensional spaces $\mathcal{X}^r_{\mathrm{h}_1} = \{v \in C^0(\bar{\Omega}_0) : v|_{\mathrm{K}} \in \mathbb{Q}_{\mathrm{r}}(K) \quad \forall K \in \mathcal{T}_{\mathrm{h}_1}\}$ and $\mathcal{X}^s_{\mathrm{h}_2} = \{v \in C^0(\bar{\Omega}_0) : v|_{\mathrm{K}} \in \mathbb{Q}_{\mathrm{s}}(K) \ \ \forall K \in \mathcal{T}_{\mathrm{h}_2}\}$, for $r, s \geq 1$.

We adopt the following notation. We denote by, e.g., $\mathbf{d}_{\mathrm{h}_1}(t) \approx \mathbf{d}(t)$ the semi-discretized Finite Element approximation of the variable $\mathbf{d}(t)$, defined over the computational mesh $\mathcal{T}_{\mathrm{h}_1}$. On the other hand, we denote by $\underline{\mathbf{d}}_{\mathrm{h}_1}(t)$ the vector collecting the DOFs associated with $\mathbf{d}_{\mathrm{h}_1}(t)$. Finally, we denote by $\underline{\mathbf{d}}^n_{\mathrm{h}_1} \simeq \underline{\mathbf{d}}_{\mathrm{h}_1}(t^n)$ the vector collecting the DOFs of the fully discretized Finite Element problem.

**Monodomain equation.** The set of basis functions for $\mathcal{X}^r_{\mathrm{h}_1}$ with $N_{\mathrm{u}} = \dim(\mathcal{X}^r_{\mathrm{h}_1})$ is given by $\{\phi_{\mathrm{i}}\}_{\mathrm{i}=1}^{N_{\mathrm{u}}}$. The semi-discretized formulation of the monodomain equation reads: find $u_{\mathrm{h}_1}(t) \in \mathcal{X}^r_{\mathrm{h}_1}$ for all $t \in (0,T)$ such that

$$\begin{aligned}
&\int_{\Omega_0} \dot{u}_{\mathrm{h}_1}(t)\phi_{\mathrm{i}}\, d\Omega_0 + \int_{\Omega_0} (J_{\mathrm{h}_1}\mathbf{F}^{-1}_{\mathrm{h}_1}(\mathbf{d}_{\mathrm{h}_1}(t))\widetilde{\boldsymbol{D}}_{\mathrm{M}}\mathbf{F}^{-T}_{\mathrm{h}_1}(\mathbf{d}_{\mathrm{h}_1}(t))\nabla u_{\mathrm{h}_1}(t)) \cdot \nabla\phi_{\mathrm{i}}\, d\Omega_0 \\
&+ \int_{\Omega_0} \widetilde{\mathcal{I}}_{\mathrm{ion}}(u_{\mathrm{h}_1}(t), \boldsymbol{w}_{\mathrm{h}_1}(t), \boldsymbol{z}_{\mathrm{h}_1}(t))\phi_{\mathrm{i}}\, d\Omega_0 = \int_{\Omega_0} \widetilde{\mathcal{I}}_{\mathrm{app}}(t)\phi_{\mathrm{i}}\, d\Omega_0 \quad \forall i = 1, ..., N_{\mathrm{u}},
\end{aligned} \tag{2}$$

with $u_{\mathrm{h}_1}(0) = \sum_{\mathrm{j}=1}^{N_{\mathrm{u}}} (u_0, \phi_{\mathrm{j}})_{L^2(\Omega_0)}\phi_{\mathrm{j}}$. The functions $\boldsymbol{w}_{\mathrm{h}_1}(t)$ and $\boldsymbol{z}_{\mathrm{h}_1}(t)$ are the semi-discretized versions of the gating variables and of the concentration variables, whereas $u_{\mathrm{h}_1}(t) = \sum_{\mathrm{j}=1}^{N_{\mathrm{u}}} u_{\mathrm{j},\mathrm{h}_1}(t)\phi_{\mathrm{j}}$ is the Finite Element solution that approximates $u = u(t)$. The tensor $\mathbf{F}_{\mathrm{h}_1}$ is the interpolated deformation tensor, obtained through the following procedure:

- The approximate solution $\mathbf{d}_{\mathrm{h}_2}$ to problem ($\mathscr{M}$) is interpolated from $\mathcal{T}_{\mathrm{h}_2}$ to $\mathcal{T}_{\mathrm{h}_1}$ to compute $\mathbf{d}_{\mathrm{h}_1}$. The evaluation of $\mathbf{d}_{\mathrm{h}_1}$ at the quadrature points used to approximate Eq. (2) is performed using our intergrid transfer operator between nested meshes and arbitrary Finite Element spaces.
- We build $\mathbf{F}_{\mathrm{h}_1} = \mathbf{I} + \nabla\mathbf{d}_{\mathrm{h}_1}$ directly on $\mathcal{T}_{\mathrm{h}_1}$ after the interpolation procedure applied on the displacement field.

We rewrite Eq. (2) as a system of nonlinear ODEs by setting $\underline{\boldsymbol{u}}_{\mathrm{h}_1}(t) = \{u_{\mathrm{j},\mathrm{h}_1}(t)\}_{\mathrm{j}=1}^{N_{\mathrm{u}}}$:

$$\begin{cases}
\mathcal{M}\dot{\underline{\boldsymbol{u}}}_{\mathrm{h}_1}(t) + \mathcal{K}(\underline{\mathbf{d}}_{\mathrm{h}_1}(t))\underline{\boldsymbol{u}}_{\mathrm{h}_1}(t) + \boldsymbol{I}_{\mathrm{ion}}(\underline{\boldsymbol{u}}_{\mathrm{h}_1}(t), \underline{\boldsymbol{w}}_{\mathrm{h}_1}(t), \underline{\boldsymbol{z}}_{\mathrm{h}_1}(t)) = \boldsymbol{I}_{\mathrm{app}}(t) \quad \forall t \in (0,T), \\
\underline{\boldsymbol{u}}_{\mathrm{h}_1}(0) = \underline{\boldsymbol{u}}_{0,\mathrm{h}_1},
\end{cases} \tag{3}$$

where we have defined the following matrices

$$\mathcal{M}_{\mathrm{ij}} = \int_{\Omega_0} \phi_{\mathrm{j}} \phi_{\mathrm{i}} \, d\Omega_0, \quad \mathcal{K}_{\mathrm{ij}}(\underline{\mathbf{d}}_{\mathrm{h}_1}(t)) = \int_{\Omega_0} (J_{\mathrm{h}_1} \mathbf{F}_{\mathrm{h}_1}^{-1} \boldsymbol{D}_m \mathbf{F}_{\mathrm{h}_1}^{-T} \nabla \phi_{\mathrm{j}}) \cdot \nabla \phi_{\mathrm{i}} \, d\Omega_0,$$

and the following vectors

$$\left(\boldsymbol{I}_{\mathrm{ion}}(\underline{\boldsymbol{u}}_{\mathrm{h}_1}(t), \underline{\boldsymbol{w}}_{\mathrm{h}_1}(t), \underline{\boldsymbol{z}}_{\mathrm{h}_1}(t))\right)_{\mathrm{i}} = \int_{\Omega_0} \widetilde{\mathcal{I}}_{\mathrm{ion}}(u_{\mathrm{h}_1}(t), \boldsymbol{w}_{\mathrm{h}_1}(t), \boldsymbol{z}_{\mathrm{h}_1}(t)) \phi_{\mathrm{i}} \, d\Omega_0,$$

$$(\boldsymbol{I}_{\mathrm{app}}(t))_{\mathrm{i}} = \int_{\Omega_0} \widetilde{\mathcal{I}}_{\mathrm{app}}(t) \phi_{\mathrm{i}} \, d\Omega_0.$$

For the evaluation of the nonlinear term $\boldsymbol{I}_{\mathrm{ion}}(\underline{\boldsymbol{u}}_{\mathrm{h}_1}(t), \underline{\boldsymbol{w}}_{\mathrm{h}_1}(t), \underline{\boldsymbol{z}}_{\mathrm{h}_1}(t))$, three strategies are available [22, 33, 36–38]. In this work, we use the so-called ionic current interpolation (ICI) approach, which yields a faster assembly of the ionic term [22, 33]. By denoting by $\{\boldsymbol{x}_q^K\}_{q=1}^{N_{\mathrm{q}}}$ and $\{\omega_q^K\}_{q=1}^{N_{\mathrm{q}}}$ the quadrature nodes and weights of a generic mesh element of $K \in \mathcal{T}_{\mathrm{h}_1}$, the term $\boldsymbol{I}_{\mathrm{ion}}(\underline{\boldsymbol{u}}_{\mathrm{h}_1}(t), \underline{\boldsymbol{w}}_{\mathrm{h}_1}(t), \underline{\boldsymbol{z}}_{\mathrm{h}_1}(t))$ is firstly evaluated at the DOFs and then interpolated at the quadrature nodes, i.e.:

$$\begin{aligned} &\int_{\Omega_0} \widetilde{\mathcal{I}}_{\mathrm{ion}}(u_{\mathrm{h}_1}(t), \boldsymbol{w}_{\mathrm{h}_1}(t), \boldsymbol{z}_{\mathrm{h}_1}(t)) \phi_{\mathrm{i}} \, d\Omega_0 \\ &\quad \approx \sum_{K \in \mathcal{T}_{\mathrm{h}_1}} \left( \sum_{q=1}^{N_q} \sum_{\mathrm{j}=1}^{N_u} \widetilde{\mathcal{I}}_{\mathrm{ion}}\left(u_{\mathrm{j},\mathrm{h}_1}(t), \boldsymbol{w}_{\mathrm{j},\mathrm{h}_1}(t) \boldsymbol{z}_{\mathrm{j},\mathrm{h}_1}(t)\right) \phi_{\mathrm{j}}(\boldsymbol{x}_q^K) \phi_{\mathrm{i}}(\boldsymbol{x}_q^K) \omega_q^K \right). \end{aligned} \tag{4}$$

**Ionic model.** The ionic model under consideration is a system of 18 ODEs (12 for the gating variables, 6 for the concentration variables), which indirectly depends on the space variable over the mesh $\mathcal{T}_{\mathrm{h}_1}$ through the transmembrane potential $u$. The semi-discrete formulation can be written as follows:

$$\begin{cases} \dot{\underline{\boldsymbol{w}}}_{\mathrm{h}_1}(t) = \overline{\boldsymbol{H}}(\underline{\boldsymbol{u}}_{\mathrm{h}_1}(t), \underline{\boldsymbol{w}}_{\mathrm{h}_1}(t)) & \forall t \in (0,T), \\ \dot{\underline{\boldsymbol{z}}}_{\mathrm{h}_1}(t) = \overline{\boldsymbol{G}}(\underline{\boldsymbol{u}}_{\mathrm{h}_1}(t), \underline{\boldsymbol{w}}_{\mathrm{h}_1}(t), \underline{\boldsymbol{z}}_{\mathrm{h}_1}(t)) & \forall t \in (0,T), \\ \underline{\boldsymbol{w}}_{\mathrm{h}_1}(0) = \underline{\boldsymbol{w}}_{0,\mathrm{h}_1}, \\ \underline{\boldsymbol{z}}_{\mathrm{h}_1}(0) = \underline{\boldsymbol{z}}_{0,\mathrm{h}_1}. \end{cases} \tag{5}$$

**Mechanical activation model.** The semi-discrete formulation, which is written on $\mathcal{T}_{\mathrm{h}_2}$, reads:

$$\begin{cases} \dot{\underline{\mathbf{s}}}_{\mathrm{h}_2}(t) = \overline{\boldsymbol{K}}(\underline{\mathbf{s}}_{\mathrm{h}_2}(t), ([\mathrm{Ca}^{2+}]_{\mathrm{i},\mathrm{h}_2}(t), SL_{\mathrm{h}_2}(t))^T) & \forall t \in (0,T), \\ \underline{\mathbf{s}}_{\mathrm{h}_2}(0) = \underline{\mathbf{s}}_{0,\mathrm{h}_2}. \end{cases} \tag{6}$$

where $[\mathrm{Ca}^{2+}]_{\mathrm{i},\mathrm{h}_2}(t)$ is obtained by interpolating the intracellular calcium concentration of TTP06 model from $\mathcal{T}_{\mathrm{h}_1}$ to $\mathcal{T}_{\mathrm{h}_2}$ as explained in the introduction to this section. The operator $\overline{\boldsymbol{K}}$ represents the element-wise application of the ANN associated with the activation model, previously trained from a collection of simulations obtained with the high-fidelity RDQ18 model [46, 48]. On the other hand, $SL_{\mathrm{h}_2}(t)$ is obtained by solving, for any $t \in (0,T)$:

$$\begin{cases} \left( SL_{\mathrm{h}_2}(t) - SL_0 \sqrt{\mathcal{I}_{4f,\mathrm{h}_2}(t)} \right) - \delta_{SL}^2 \Delta SL_{\mathrm{h}_2}(t) = 0 & \text{in } \Omega_0 \times (0,T) \\ \delta_{SL}^2 \nabla SL_{\mathrm{h}_2}(t) \cdot \mathbf{N}_{\mathrm{h}_2} = 0 & \text{on } \partial\Omega_0 \times (0,T) \end{cases} \tag{7}$$

where $\mathcal{I}_{4f,\mathrm{h}_2}(t) = \mathbf{F}_{\mathrm{h}_2}(t)\mathbf{f}_0 \cdot \mathbf{F}_{\mathrm{h}_2}(t)\mathbf{f}_0$. Finally, $T_{\mathrm{a},\mathrm{h}_2}(t)$ denotes the semi-discretized active tension, obtained by evaluating the function $T_{\mathrm{a}} = T_{\mathrm{a}}^{\max} G(\mathbf{s})$ at the nodes.

**Mechanical model.** We denote by $[\mathcal{X}_{h_2}^s]^3$ the finite dimensional subspace of vector valued functions and by $\{\boldsymbol{\phi}_i\}_{i=1}^{N_\mathbf{d}}$ its basis. The semi-discretized version of ($\mathscr{M}$) reads: given $T_{a,h_2}(t)$, find $\mathbf{d}_{h_2} = \mathbf{d}_{h_2}(t) \in [\mathcal{X}_{h_2}^s]^3$ for all $t \in (0, T)$ such that

$$\begin{aligned}
&\int_{\Omega_0} \rho_s \ddot{\mathbf{d}}_{h_2}(t) \cdot \boldsymbol{\phi}_i \, d\Omega_0 + \int_{\Omega_0} \mathbf{P}(\mathbf{d}_{h_2}(t), T_{a,h_2}(t)) : \nabla \boldsymbol{\phi}_i d\Omega_0 \\
&\quad + \int_{\Gamma_0^{\text{epi}}} \Big[ (\mathbf{N}_{h_2} \otimes \mathbf{N}_{h_2}) \Big( K_\perp^{\text{epi}} \mathbf{d}_{h_2}(t) + C_\perp^{\text{epi}} \dot{\mathbf{d}}_{h_2}(t) \Big) \\
&\quad + (\mathbf{I} - \mathbf{N}_{h_2} \otimes \mathbf{N}_{h_2}) \Big( K_\parallel^{\text{epi}} \mathbf{d}_{h_2}(t) + C_\parallel^{\text{epi}} \dot{\mathbf{d}}_{h_2}(t) \Big) \Big] \cdot \boldsymbol{\phi}_i \, d\Gamma_0 \\
&\quad = -p_{\text{LV}}(t) \int_{\Gamma_0^{\text{endo}}} J_{h_2} \mathbf{F}_{h_2}^{-T} \mathbf{N} \cdot \boldsymbol{\phi}_i \, d\Gamma_0 + p_{\text{LV}}(t) \int_{\Gamma_0^{\text{base}}} |J_{h_2} \mathbf{F}_{h_2}^{-T} \mathbf{N}| \mathbf{v}_{h_2}^{\text{base}} \cdot \boldsymbol{\phi}_i \, d\Gamma_0 \\
&\quad \forall i = 1, ..., N_\mathbf{d},
\end{aligned} \tag{8}$$

with $\mathbf{d}_{h_2}(0) = \sum_{j=1}^{N_\mathbf{d}} (\mathbf{d}_0, \boldsymbol{\phi}_j)_{[L^2(\Omega_0)]^3} \boldsymbol{\phi}_j$, $\dot{\mathbf{d}}_{h_2}(0) = \sum_{j=1}^{N_\mathbf{d}} (\dot{\mathbf{d}}_0, \boldsymbol{\phi}_j)_{[L^2(\Omega_0)]^3} \boldsymbol{\phi}_j$ and where

$$\mathbf{v}_{h_2}^{\text{base}} = \frac{\int_{\Gamma_0^{\text{endo}}} J_{h_2} \mathbf{F}_{h_2}^{-T} \mathbf{N} \, d\Gamma_0}{\int_{\Gamma_0^{\text{base}}} |J_{h_2} \mathbf{F}_{h_2}^{-T} \mathbf{N}| \, d\Gamma_0}.$$

The corresponding algebraic formulation reads:

$$\begin{cases}
\rho_s \mathcal{M} \ddot{\underline{\mathbf{d}}}_{h_2}(t) + \mathcal{F} \dot{\underline{\mathbf{d}}}_{h_2}(t) + \mathcal{G} \underline{\mathbf{d}}_{h_2}(t) + \boldsymbol{S}(\underline{\mathbf{d}}_{h_2}(t), T_{a,h_2}(t)) = p_{\text{LV}}(t) \, \boldsymbol{p}(\underline{\mathbf{d}}_{h_2}(t)) \\
\underline{\mathbf{d}}_{h_2}(0) = \underline{\mathbf{d}}_{0,h_2}, \quad \dot{\underline{\mathbf{d}}}_{h_2}(0) = \dot{\underline{\mathbf{d}}}_{0,h_2},
\end{cases} \quad \forall t \in (0, T), \tag{9}$$

with:

$$\boldsymbol{S}_i(\underline{\mathbf{d}}_{h_2}(t), T_{a,h_2}(t)) = \int_{\Omega_0} \mathbf{P}(\mathbf{d}_{h_2}(t), T_{a,h_2}(t)) : \nabla \boldsymbol{\phi}_i \, d\Omega_0,$$

$$\mathcal{F}_{i,j} = \int_{\Gamma_0^{\text{epi}}} \Big[ (\mathbf{N}_{h_2} \otimes \mathbf{N}_{h_2}) C_\perp^{\text{epi}} + (\mathbf{I} - \mathbf{N}_{h_2} \otimes \mathbf{N}_{h_2}) C_\parallel^{\text{epi}} \Big] \boldsymbol{\phi}_j \cdot \boldsymbol{\phi}_i \, d\Gamma_0,$$

$$\mathcal{G}_{i,j} = \int_{\Gamma_0^{\text{epi}}} \Big[ (\mathbf{N}_{h_2} \otimes \mathbf{N}_{h_2}) K_\perp^{\text{epi}} + (\mathbf{I} - \mathbf{N}_{h_2} \otimes \mathbf{N}_{h_2}) K_\parallel^{\text{epi}} \Big] \boldsymbol{\phi}_j \cdot \boldsymbol{\phi}_i \, d\Gamma_0,$$

$$\boldsymbol{p}_i(\underline{\mathbf{d}}_{h_2}(t)) = \int_{\Gamma_0^{\text{base}}} |J_{h_2} \mathbf{F}_{h_2}^{-T} \mathbf{N}| \mathbf{v}_{h_2}^{\text{base}} \cdot \boldsymbol{\phi}_i \, d\Gamma_0 - \int_{\Gamma_0^{\text{endo}}} J_{h_2} \mathbf{F}_{h_2}^{-T} \mathbf{N} \cdot \boldsymbol{\phi}_i \, d\Gamma_0,$$

where $\underline{\mathbf{d}}_{h_2}(t) = \{\mathbf{d}_{j,h_2}(t)\}_{j=1}^{N_\mathbf{d}}$. Equations (2), (5), (6) and (8) provide a splitted semi-discretization of the entire electromechanical model.

## 3.2 Time discretization

We define $\Delta t$ as the time step for ($\mathscr{M}_I$)–($\mathscr{V}_I$) and $\tau = \Delta t / N_{\text{sub}}$ as the time step for ($\mathscr{I}_{SI}$), ($\mathscr{E}_{SI}$), ($\mathscr{A}_E$) and ($\mathscr{C}_E$), being $N_{\text{sub}} \in \mathbb{N}$ the number of intermediate substeps that need to be solved by ($\mathscr{I}_{SI}$)–($\mathscr{E}_{SI}$)–($\mathscr{A}_E$)–($\mathscr{C}_E$) before a time step $\Delta t$ of ($\mathscr{M}_I$)–($\mathscr{V}_I$) is performed.

Problem ($\mathscr{I}_{SI}$)–($\mathscr{E}_{SI}$)–($\mathscr{A}_E$)–($\mathscr{C}_E$) from $t^n$ to $t^{n+1}$, once we set $t^{n+\frac{m}{N_{\text{sub}}}} = t^n + m\tau$, for $m = 1, ..., N_{\text{sub}}$, reads as follows:

- We find $\underline{\boldsymbol{w}}_{\mathrm{h}_1}^{n+\frac{m}{N_{\mathrm{sub}}}}$ and $\underline{\boldsymbol{z}}_{\mathrm{h}_1}^{n+\frac{m}{N_{\mathrm{sub}}}}$ defined on $\mathcal{T}_{\mathrm{h}_1}$ by solving:

$$\begin{cases} \frac{1}{\tau}\underline{\boldsymbol{w}}_{\mathrm{h}_1}^{n+\frac{m}{N_{\mathrm{sub}}}} = \frac{1}{\tau}\underline{\boldsymbol{w}}_{\mathrm{h}_1}^{n} + \overline{\boldsymbol{H}}(\underline{\boldsymbol{u}}_{\mathrm{h}_1}^{n}, \underline{\boldsymbol{w}}_{\mathrm{h}_1}^{n+\frac{m}{N_{\mathrm{sub}}}}), \\ \frac{1}{\tau}\underline{\boldsymbol{z}}_{\mathrm{h}_1}^{n+\frac{m}{N_{\mathrm{sub}}}} = \frac{1}{\tau}\underline{\boldsymbol{z}}_{\mathrm{h}_1}^{n} + \overline{\boldsymbol{G}}(\underline{\boldsymbol{u}}_{\mathrm{h}_1}^{n}, \underline{\boldsymbol{w}}_{\mathrm{h}_1}^{n}, \underline{\boldsymbol{z}}_{\mathrm{h}_1}^{n}). \end{cases} \tag{10}$$

  We adopt here the first order implicit-explicit (IMEX) scheme proposed in [45]. Specifically, we employ an explicit treatment of the ionic concentrations to avoid the solution of a nonlinear system (such choice does not compromise the stability of the scheme, thanks to the non-stiff dynamics of concentrations), and an implicit treatment of the gating variables, because of the severe CFL condition on the time step induced by an explicit scheme. We notice that, thanks to the linear dynamics of the gating variables, such implicit handling does not require the solution of a system of linear or nonlinear equations.

- We interpolate $\underline{\mathbf{d}}_{\mathrm{h}_2}^{n}$ on the fine mesh $\mathcal{T}_{\mathrm{h}_1}$ once per time step, at $t = t^n$. We use $\underline{\boldsymbol{z}}_{\mathrm{h}_1}^{n+\frac{m}{N_{\mathrm{sub}}}}$ from (10) and $\underline{\mathbf{d}}_{\mathrm{h}_1}^{n}$ to find $\underline{\boldsymbol{u}}_{\mathrm{h}_1}^{n+\frac{m}{N_{\mathrm{sub}}}}$ over $\mathcal{T}_{\mathrm{h}_1}$ by solving:

$$\begin{aligned} &\left(\frac{1}{\tau}\mathcal{M} + \mathcal{K}(\underline{\mathbf{d}}_{\mathrm{h}_1}^{n}) + \boldsymbol{I}_{\mathrm{u}}^{\mathrm{ion}}\left(\underline{\boldsymbol{u}}_{\mathrm{h}_1}^{n}, \underline{\boldsymbol{z}}_{\mathrm{h}_1}^{n+\frac{m}{N_{\mathrm{sub}}}}\right)\right)\underline{\boldsymbol{u}}_{\mathrm{h}_1}^{n+\frac{m}{N_{\mathrm{sub}}}} = \\ &\qquad \frac{1}{\tau}\mathcal{M}\underline{\boldsymbol{u}}_{\mathrm{h}_1}^{n} - \tilde{\boldsymbol{I}}^{\mathrm{ion}}\left(\underline{\boldsymbol{u}}_{\mathrm{h}_1}^{n}, \underline{\boldsymbol{z}}_{\mathrm{h}_1}^{n+\frac{m}{N_{\mathrm{sub}}}}\right) + \boldsymbol{I}_{\mathrm{app}}\left(t^{n+\frac{m}{N_{\mathrm{sub}}}}\right). \end{aligned} \tag{11}$$

  $\boldsymbol{I}_{\mathrm{u}}^{\mathrm{ion}}$ is the derivative of the terms of $\boldsymbol{I}_{\mathrm{ion}}$ that linearly depends on $\underline{\boldsymbol{u}}_{\mathrm{h}_1}$, while $\tilde{\boldsymbol{I}}^{\mathrm{ion}}$ collects all the other terms.

- We interpolate $[\mathrm{Ca}^{2+}]_{\mathrm{i},\mathrm{h}_2}^{n+\frac{m}{N_{\mathrm{sub}}}}$ from (10) on the coarse mesh $\mathcal{T}_{\mathrm{h}_2}$, and we find $\underline{\mathbf{s}}_{\mathrm{h}_2}^{n+\frac{m}{N_{\mathrm{sub}}}}$ by solving:

$$\underline{\mathbf{s}}_{\mathrm{h}_2}^{n+\frac{m}{N_{\mathrm{sub}}}} = \underline{\mathbf{s}}_{\mathrm{h}_2}^{n} + \tau\overline{\boldsymbol{K}}(\underline{\mathbf{s}}_{\mathrm{h}_2}^{n}, ([\mathrm{Ca}^{2+}]_{\mathrm{i},\mathrm{h}_2}^{n+\frac{m}{N_{\mathrm{sub}}}}, SL_{\mathrm{h}_2}^{n})^T). \tag{12}$$

  where $SL_{\mathrm{h}_2}^{n}$ is obtained by solving problem (7).

- We find $\boldsymbol{c}_1^{n+\frac{m}{N_{\mathrm{sub}}}}$, with the $4^{th}$ order Runge-Kutta method:

$$\boldsymbol{c}_1^{n+\frac{m}{N_{\mathrm{sub}}}} = \boldsymbol{c}_1^{n} + \frac{1}{6}\left(\mathbf{k}_1 + 2\mathbf{k}_2 + 2\mathbf{k}_3 + \mathbf{k}_4\right), \tag{13}$$

  with $\mathbf{k}_1$, $\mathbf{k}_2$, $\mathbf{k}_3$ and $\mathbf{k}_4$ computed as follows:

$$\begin{aligned} &\mathbf{k}_1 = \tau\widetilde{\boldsymbol{D}}\left(t^n, \boldsymbol{c}_1^n, p_{\mathrm{LV}}^n\right), &&\mathbf{k}_2 = \tau\widetilde{\boldsymbol{D}}\left(t^{n+\frac{1}{2}\frac{m}{N_{\mathrm{sub}}}}, \boldsymbol{c}_1^n + \frac{\mathbf{k}_1}{2}, p_{\mathrm{LV}}^n\right), \\ &\mathbf{k}_3 = \tau\widetilde{\boldsymbol{D}}\left(t^{n+\frac{1}{2}\frac{m}{N_{\mathrm{sub}}}}, \boldsymbol{c}_1^n + \frac{\mathbf{k}_2}{2}, p_{\mathrm{LV}}^n\right), &&\mathbf{k}_4 = \tau\widetilde{\boldsymbol{D}}\left(t^{n+\frac{m}{N_{\mathrm{sub}}}}, \boldsymbol{c}_1^n + \mathbf{k}_3, p_{\mathrm{LV}}^n\right). \end{aligned}$$

After having solved (10), (11), (12) and (13) for $N_{\mathrm{sub}}$ steps, we treat $(\mathcal{M}_I)$–$(\mathscr{V}_{\mathrm{I}})$ at $t^{n+1}$ by updating $\underline{\mathbf{d}}_{\mathrm{h}_2}^{n+1}$ and $p_{\mathrm{LV}}^{n+1}$ with the following system:

$$\begin{cases} \left(\rho_{\mathrm{s}}\frac{1}{\Delta t^2}\mathcal{M} + \frac{1}{\Delta t}\mathcal{F} + \mathcal{G}\right)\underline{\mathbf{d}}_{\mathrm{h}_2}^{n+1} + \boldsymbol{S}(\underline{\mathbf{d}}_{\mathrm{h}_2}^{n+1}, \mathbf{T}_{\mathrm{a},\mathrm{h}_2}^{n+1}) \\ \qquad = \rho_{\mathrm{s}}\frac{2}{\Delta t^2}\mathcal{M}\underline{\mathbf{d}}_{\mathrm{h}_2}^{n} - \rho_{\mathrm{s}}\frac{1}{\Delta t^2}\mathcal{M}\underline{\mathbf{d}}_{\mathrm{h}_2}^{n-1} + \frac{1}{\Delta t}\mathcal{F}\underline{\mathbf{d}}_{\mathrm{h}_2}^{n} + p_{\mathrm{LV}}^{n+1}\boldsymbol{p}(\underline{\mathbf{d}}_{\mathrm{h}_2}^{n}, \underline{\mathbf{d}}_{\mathrm{h}_2}^{n+1}), \\ V_{\mathrm{LV}}^{\mathrm{3D}}(\underline{\mathbf{d}}_{\mathrm{h}_2}^{n+1}) = V_{\mathrm{LV}}^{\mathrm{0D}}(\boldsymbol{c}_1^{n+1}). \end{cases} \tag{14}$$

Equation (14) is a nonlinear saddle-point problem. In the next section (Sec. 3.3) we provide details about its numerical approximation at the algebraic level. An alternative to our approach for numerical discretization would be employing implicit monolithic schemes, which are known to be stable and accurate, but at the same time they present constraints in the choice of the time steps too and they are characterized by high computational costs [22, 52]. Indeed, we can only use one time step and one mesh with a monolithic approach. Therefore, we are forced to choose a small time step and a fine representation of the computational domain due to the requirements of cardiac electrophysiology. Our approach is instead accurate and computationally efficient, thanks to the flexibility in the choice of both space and time resolution among the different core models.

### 3.3 Algorithm for the resolution of Eq. (14)

We approximate the solution of Eq. (14) by means of a quasi-Newton strategy [42], as we proposed in [48]. Specifically, in the computation of the Jacobian matrix, we neglect the derivative of the nonlocal term $\mathbf{v}_{\mathrm{h}_2}^{\mathrm{base}}$ in the pressure variable, and we update the Jacobian at each time step, but not through the iterations of Newton's loop. By moving all the terms in Eq. (14) to the left hand side and by rewriting its first and second line as $\mathbf{r}_{\mathbf{d}}^{n+1}(\underline{\mathbf{d}}_{\mathrm{h}_2}^{n+1}, p_{\mathrm{LV}}^{n+1}) = \mathbf{0}$ and $r_{\mathrm{p}}^{n+1}(\underline{\mathbf{d}}_{\mathrm{h}_2}^{n+1}) = 0$, respectively, the quasi-Newton algorithm reads as follows:

- We set $\underline{\mathbf{d}}_{\mathrm{h}_2}^{n+1,0} = \underline{\mathbf{d}}_{\mathrm{h}_2}^{n}$ and $p_{\mathrm{LV}}^{n+1,0} = p_{\mathrm{LV}}^{n}$.
- For $j = 0, 1, \ldots$, until a convergence criterion is not fulfilled, we solve the following linear system:
$$\begin{pmatrix} J_{\mathbf{d},\mathbf{d}}^{n+1} & J_{\mathbf{d},\mathrm{p}}^{n+1} \\ J_{\mathrm{p},\mathbf{d}}^{n+1} & 0 \end{pmatrix} \begin{pmatrix} \Delta \underline{\mathbf{d}}_{\mathrm{h}_2}^{n+1,j} \\ \Delta p_{\mathrm{LV}}^{n+1,j} \end{pmatrix} = - \begin{pmatrix} \mathbf{r}_{\mathbf{d}}^{n+1,j} \\ r_{\mathrm{p}}^{n+1,j} \end{pmatrix}, \tag{15}$$
where $J_{\mathbf{d},\mathbf{d}}^{n+1} \simeq \frac{\partial}{\partial \underline{\mathbf{d}}} \mathbf{r}_{\mathbf{d}}^{n+1}(\underline{\mathbf{d}}_{\mathrm{h}_2}^{n}, p_{\mathrm{LV}}^{n})$ (wherein we neglected the derivative with respect to the nonlocal term $\mathbf{v}_{\mathrm{h}_2}^{\mathrm{base}}$), $J_{\mathbf{d},\mathrm{p}}^{n+1} = \frac{\partial}{\partial p} \mathbf{r}_{\mathbf{d}}^{n+1}(\underline{\mathbf{d}}_{\mathrm{h}_2}^{n}, p_{\mathrm{LV}}^{n})$, $J_{\mathrm{p},\mathbf{d}}^{n+1} = \frac{\partial}{\partial \underline{\mathbf{d}}} r_{\mathrm{p}}^{n+1}(\underline{\mathbf{d}}_{\mathrm{h}_2}^{n})$, $\mathbf{r}_{\mathbf{d}}^{n+1,j} = \mathbf{r}_{\mathbf{d}}^{n+1}(\underline{\mathbf{d}}_{\mathrm{h}_2}^{n+1,j}, p_{\mathrm{LV}}^{n+1,j})$ and $r_{\mathrm{p}}^{n+1,j} = r_{\mathrm{p}}^{n+1}(\underline{\mathbf{d}}_{\mathrm{h}_2}^{n+1,j})$ .
- We update $\underline{\mathbf{d}}_{\mathrm{h}_2}^{n+1,j+1} = \underline{\mathbf{d}}_{\mathrm{h}_2}^{n+1,j} + \Delta \underline{\mathbf{d}}_{\mathrm{h}_2}^{n+1,j}$ and $p_{\mathrm{LV}}^{n+1,j+1} = p_{\mathrm{LV}}^{n+1,j} + \Delta p_{\mathrm{LV}}^{n+1,j}$.
- When the convergence criterion is satisfied, we set $\underline{\mathbf{d}}_{\mathrm{h}_2}^{n+1} = \underline{\mathbf{d}}_{\mathrm{h}_2}^{n+1,j}$ and $p_{\mathrm{LV}}^{n+1} = p_{\mathrm{LV}}^{n+1,j}$.

From the algebraic viewpoint, we solve the saddle-point problem (15) via Schur complement reduction [7]. Specifically, we solve the two linear systems

$$J_{\mathbf{d},\mathbf{d}}^{n+1} \mathbf{v}^{n+1,j} = \mathbf{r}_{\mathbf{d}}^{n+1,j}, \qquad J_{\mathbf{d},\mathbf{d}}^{n+1} \mathbf{w}^{n+1,j} = J_{\mathbf{d},\mathrm{p}}^{n+1}$$

and we set:

$$\Delta p_{\mathrm{LV}}^{n+1,j} = \frac{r_p^{n+1,j} - J_{\mathrm{p},\mathbf{d}}^{n+1} \mathbf{v}^{n+1,j}}{J_{\mathrm{p},\mathbf{d}}^{n+1} \mathbf{w}^{n+1,j}}, \quad \Delta \underline{\mathbf{d}}_{\mathrm{h}_2}^{n+1,j} = - \left( \mathbf{v}^{n+1,j} + \mathbf{w}^{n+1,j} \Delta p_{\mathrm{LV}}^{n+1,j} \right). \tag{16}$$

We remark that, thanks to the reuse of the Jacobian matrix throughout the Newton loop, $\mathbf{w}^{n+1,j}$ becomes independent of $j$ and thus it does not need to be recomputed at each iteration. Our scheme only involves the following operations: for each time step, we assemble the matrix $J_{\mathbf{d},\mathbf{d}}^{n+1}$ and the vectors $J_{\mathrm{p},\mathbf{d}}^{n+1}$ and $J_{\mathbf{d},\mathrm{p}}^{n+1}$ and we solve the linear system $J_{\mathbf{d},\mathbf{d}}^{n+1} \mathbf{w}^{n+1,0} = J_{\mathbf{d},\mathrm{p}}^{n+1}$; at each Newton iteration, we only need to solve the linear system $J_{\mathbf{d},\mathbf{d}}^{n+1} \mathbf{v}^{n+1,j} = \mathbf{r}_{\mathbf{d}}^{n+1,j}$ and perform a couple of matrix-vector multiplications and a vector-vector sum of Eq. (16). We later show, through several numerical simulations, that this approach is numerically stable. Even more importantly, our scheme is appropriate for the whole heartbeat, as it does not require adaptations according to the specific cardiac phase.

Moreover, our approach allows for a segregated solution of the 3D-0D coupled model, which is instead typically solved through a monolithic strategy [28]. As a matter of fact, segregated schemes available in literature that couple a mechanical problem with a model describing the dynamics of a fluid (even when described through a 0D circulation model) generally fail when the fluid domain is fully enclosed by the solid structure, because the incompressibility constraint of the fluid is no longer satisfied after the structure update [28]. This issue, known as *balloon dilemma*, affects also cardiac chambers, either when coupled with a 3D circulation model – within a fluid-structure interaction framework [5, 25] – or with a 0D one. To overcome this issue, the 3D-0D cardiac circulation models available in literature rely either on a monolithic solution of the two models, where the 0D and 3D models are simultaneously discretized as a unique system and then typically solved by a Newton method [28] or iterative methods that progressively update the cavity pressures and the solid displacement, until convergence is reached. For instance, in [32], the authors proposed a method where the cavity pressure is initially estimated by extrapolating from previous time steps. Then, the cavity compliance (i.e. $\partial V/\partial p$, where $p$ and $V$ are the cavity pressure and volume, respectively) is estimated by finite differences and is used to update the pressure until the blood flux of the 3D model matches that of the 0D model within a prescribed tolerance. Similarly, in [17, 21], during the isovolumic phases, the cavity pressure is iteratively updated by a fixed point scheme. However, convergence of this scheme depends upon a relaxation parameter, whose optimal value needs to be manually assessed from case by case [21].

With our approach, instead, the mechanical 3D model ($\mathscr{M}$) and the circulation 0D model ($\mathscr{C}$) are solved in both a segregated and staggered manner. Indeed, we do not solve the ($\mathscr{M}$) model simultaneously to the ($\mathscr{C}$) model, but coupled to the volume-consistency condition ($\mathscr{V}$) instead. In this way we end up with the saddle-point problem (14). Thanks to the above algorithm, finally, at each Newton iteration, we only need to solve a linear system involving the Jacobian matrix of the ($\mathscr{M}$) problem, and perform a couple of matrix-vector operations.

### 3.4 Stresses computation

Starting from [21], we define the following indicator for evaluating the components of the mechanical stress:

$$S_{\mathrm{ab}} = (\mathbf{P}\mathbf{a}_0) \cdot \frac{\mathbf{F}\mathbf{b}_0}{|\mathbf{F}\mathbf{b}_0|}, \tag{17}$$

where $a, b \in \{f, s, n\}$ indicate the fibers ($f$), the sheets ($s$) and the crossfibers ($n$) direction respectively. The metric $S_{\mathrm{ab}}$ measures the stress component in the $\mathbf{b}$ direction (where $\mathbf{b} = \frac{\mathbf{F}\mathbf{b}_0}{|\mathbf{F}\mathbf{b}_0|}$ denotes the direction $\mathbf{b}_0$ in the current configuration) across a surface normal to the direction $\mathbf{a}_0$. Hence, we refer to axial stresses when $a = b$, whereas we get shear stresses when $a \neq b$. We remark that in the active stress framework, $\mathbf{P}$ incorporates an additive decomposition between the passive and the active terms, the latter coming from the active tension $T_{\mathrm{a}}$.

Even though a possible strategy to handle the calculation of the stress tensor would be to solve a $L^2$-projection problem [21], we represent each component as a piece-wise constant ($\mathbb{Q}_0$) Finite Element vector, where the average of the values over the quadrature points of each cell is associated with the only local degree of freedom corresponding to the cell centroid. Since the stresses are only processed for visualization purposes, our strategies turns out to be very efficient while not hampering the accuracy of the computation.

# 4 Recovering the reference stress-free configuration

Here we present the algorithm for recovering the reference configuration $\Omega_0$ from the deformed configuration $\widetilde{\Omega}$, knowing that the latter is obtained from the former by applying a pressure $\widetilde{p}$ and an active tension $\widetilde{T}_\mathrm{a}$. The steady state version of the PDE for cardiac mechanics reads:

$$\begin{cases} \nabla \cdot \mathbf{P}(\mathbf{d}, T_\mathrm{a}) = \mathbf{0} & \text{in } \Omega_0, \\ \mathbf{P}(\mathbf{d}, T_\mathrm{a})\mathbf{N} + \mathbf{K}^{\mathrm{epi}}\mathbf{d} = \mathbf{0} & \text{on } \Gamma_0^{\mathrm{epi}}, \\ \mathbf{P}(\mathbf{d}, T_\mathrm{a})\mathbf{N} = p_\mathrm{LV}\, |J\mathbf{F}^{-T}\mathbf{N}|\, \mathbf{v}^{\mathrm{base}}(t) & \text{on } \Gamma_0^{\mathrm{base}}, \\ \mathbf{P}(\mathbf{d}, T_\mathrm{a})\mathbf{N} = -p_\mathrm{LV}\, J\mathbf{F}^{-T}\mathbf{N} & \text{on } \Gamma_0^{\mathrm{endo}}. \end{cases} \tag{18}$$

In what follows, we denote by $\mathbf{d} = \mathbf{d}_\mathrm{eq}(\mathbf{x}_0, p_\mathrm{LV}, T_\mathrm{a})$ the equilibrium solution of Eq. (18) obtained on the computational domain of coordinate $\mathbf{x}_0$. Hence, our aim is finding a coordinate $\mathbf{x}_0$ such that $\mathbf{x}_0 + \mathbf{d}_\mathrm{eq}(\mathbf{x}_0, p_\mathrm{LV}, T_\mathrm{a}) = \widetilde{\mathbf{x}}$.

## 4.1 Algorithm for the recovery of the reference configuration

A representation of this algorithm is shown in Fig. 2. We start by setting the coordinate of the reference configuration $\mathbf{x}_0$ equal to the coordinate of the deformed one (i.e. $\mathbf{x}_0^{(0)} = \widetilde{\mathbf{x}}$). Then, we solve the elastostatic problem of Eq. (18), and we get the displacement $\mathbf{d}^{(0)} = \mathbf{d}_\mathrm{eq}(\mathbf{x}_0^{(0)}, p_\mathrm{LV}, T_\mathrm{a})$ (Fig. 2, top-left). Since, when the configuration $\widetilde{\Omega}$ is recorded the active tension $\widetilde{T}_\mathrm{a}$ is almost zero, in the deformed configuration $\mathbf{x}^{(0)} = \mathbf{x}_0^{(0)} + \mathbf{d}^{(0)}$ the ventricle is inflated compared to the configuration $\widetilde{\mathbf{x}}$. Thus, with the aim of correcting the mismatch between $\mathbf{x}^{(0)}$ and $\widetilde{\mathbf{x}}$, we deflate the ventricle by setting $\mathbf{x}_0^{(1)} = \mathbf{x}_0^{(0)} + (\widetilde{\mathbf{x}} - \mathbf{x}^{(0)}) = \widetilde{\mathbf{x}} - \mathbf{d}^{(0)}$ (Fig. 2, top-right). Then, we proceed by iterating the above steps. More precisely, for $k \geq 1$, we compute $\mathbf{d}^{(k)} = \mathbf{d}_\mathrm{eq}(\mathbf{x}_0^{(k)}, p_\mathrm{LV}, T_\mathrm{a})$ (Fig. 2, bottom-left) and we set $\mathbf{x}_0^{(k+1)} = \widetilde{\mathbf{x}} - \mathbf{d}^{(k)}$ (Fig. 2, bottom-right), stopping when the difference between two consecutive iterations is sufficiently small.

The whole procedure is reported in Algorithm 1, where we denote the function that solves problem (18) on the geometry with coordinates $\mathbf{x}_0$ as SteadyStateMechanics. The latter function solves the nonlinear system by means of the Newton method. In case the Newton iterations do not reach convergence (according to a criterion based both on the residual and on the difference between consecutive iterations), it returns a flag to indicate failure of the algorithm. More precisely, the function signature reads

$$(\texttt{converged_SSM}, \mathbf{d}) = \textsc{SteadyStateMechanics}(\mathbf{x}_0, p_\mathrm{LV}, T_\mathrm{a})$$

where, in case of convergence, `converged_SSM` is true and $\mathbf{d} = \mathbf{d}_\mathrm{eq}(\mathbf{x}_0, p_\mathrm{LV}, T_\mathrm{a})$, while in case of non convergence `converged_SSM` is false and $\mathbf{d}$ is not used. We remark that an algorithm similar to Algorithm 1 is presented in [43].

Algorithm 1 can be interpreted as a fixed-point iteration scheme. Indeed, the fixed-point $\mathbf{x}_0$ of the map defined by the Algorithm 1 iteration satisfies $\mathbf{x}_0 = \widetilde{\mathbf{x}} - \mathbf{d}_\mathrm{eq}(\mathbf{x}_0, p_\mathrm{LV}, T_\mathrm{a})$. The fixed-point iterations of Algorithm 1 are sketched in Fig. 3. The solution is obtained as the intersection between the line $\widetilde{\mathbf{x}} = \mathbf{x}_0 + \mathbf{d}$ and the manifold $\mathbf{d} = \mathbf{d}_\mathrm{eq}(\mathbf{x}_0, p_\mathrm{LV}, T_\mathrm{a})$. The algorithm proceeds iteratively in the space $(\mathbf{x}_0, \mathbf{d})$: the first variable is updated by the fixed-point iterations (horizontal axis of Fig. 3), while the second variable is updated by Newton iterations (vertical axis of Fig. 3).

However, the implementation shown in Algorithm 1 has several limitations when applied to realistic heart geometries and to highly nonlinear constitutive laws, such as in the case of cardiac

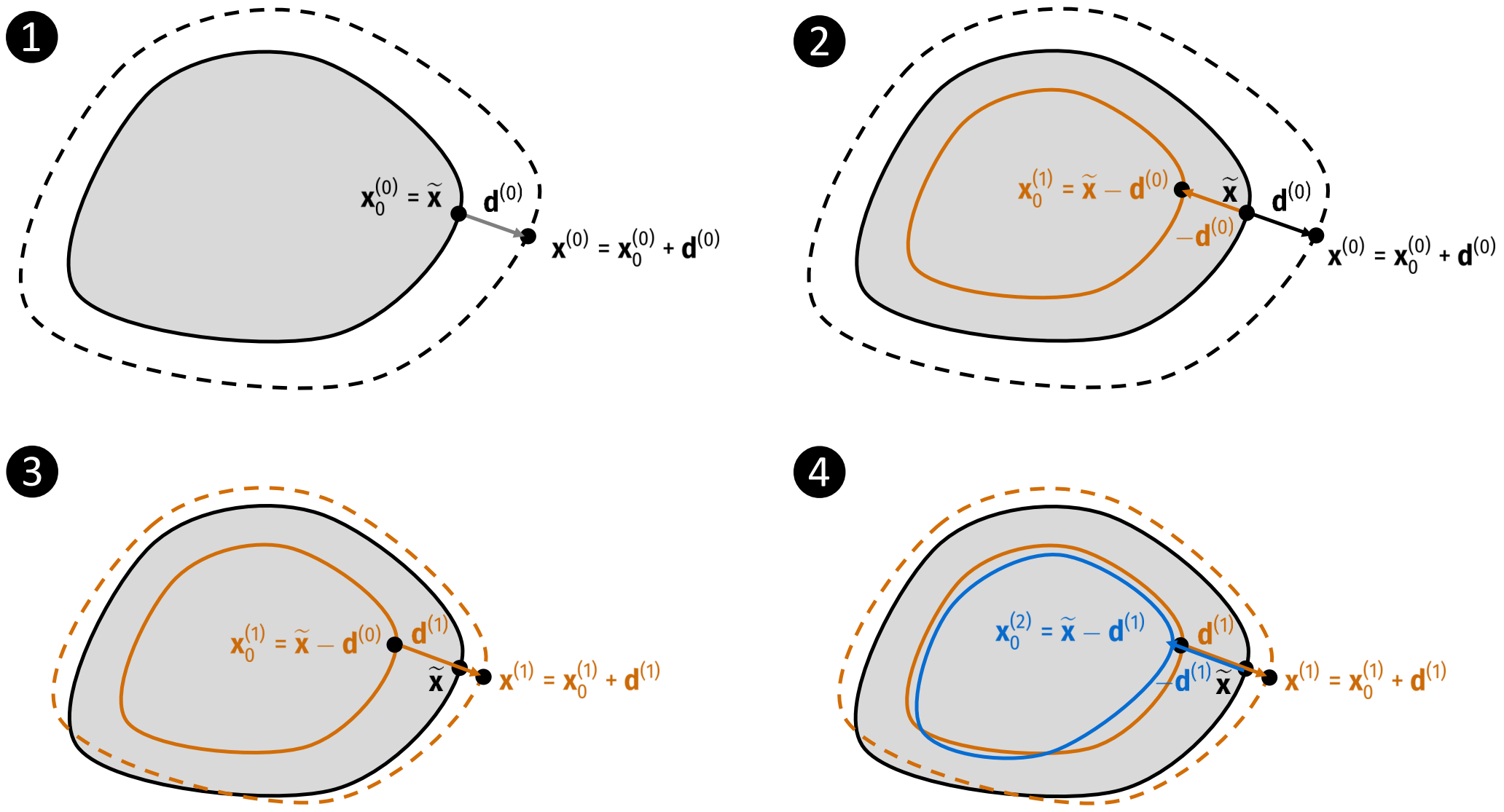

Figure 2: Visual representation of the basic version of the fixed-point algorithm (see Algorithm 1).

---

**Algorithm 1** Reference configuration recovery (basic version)

---

**parameters:** $k_{\max}, \epsilon_{\mathrm{tol}}$
**output:** `converged_RCB`, $\mathbf{x}_0$

**procedure** ReferenceConfigurationBase($\widetilde{\mathbf{x}}, \widetilde{p}, \widetilde{T}_{\mathrm{a}}$)
  $\mathbf{x}_0^{(0)} \leftarrow \widetilde{\mathbf{x}}$
  **for** $k = 0, \ldots, k_{\max}$ **do**
    (`converged_SSM`, $\mathbf{d}^{(k)}$) $\leftarrow$ SteadyStateMechanics($\mathbf{x}_0^{(k)}, \widetilde{p}, \widetilde{T}_{\mathrm{a}}$)
    **if not** `converged_SSM` **then**
      **return** (false, $\mathbf{0}$) ▷ Newton method does not converge.
    **end if**
    $\mathbf{x}^{(k)} \leftarrow \mathbf{x}_0^{(k)} + \mathbf{d}^{(k)}$
    **if** $\|\mathbf{x}^{(k)} - \widetilde{\mathbf{x}}\| \leq \epsilon_{\mathrm{tol}} \|\mathbf{d}^{(k)}\|$ **then**
      **return** (true, $\mathbf{x}_0^{(k)}$) ▷ Fixed-point converged.
    **end if**
    $\mathbf{x}_0^{(k+1)} \leftarrow \widetilde{\mathbf{x}} - \mathbf{d}^{(k)}$ ▷ Fixed-point update.
  **end for**
  **return** (false, $\mathbf{0}$) ▷ Maximum number of iterations reached.
**end procedure**

---

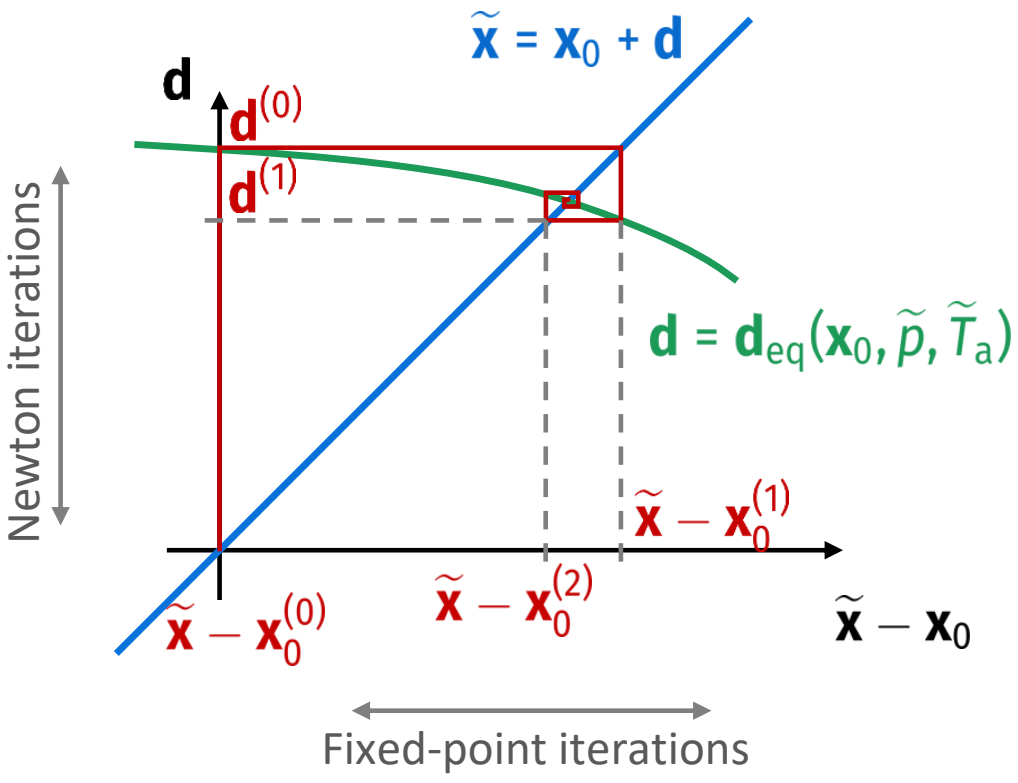

Figure 3: Representation of the basic version of the reference configuration recovery algorithm (Algorithm 1).

mechanics. In fact, both the Newton method employed in SteadyStateMechanics and the fixed-point scheme employed in ReferenceConfigurationBase do not converge if the initial guess is not “sufficiently” close to the solution. In particular, our experience revealed that the attraction basin of the fixed-point iterations scheme gets smaller and smaller as the value of $p_{\mathrm{LV}}$ grows, which makes the base formulation of Algorithm 1 not suitable for realistic cardiac geometries.

To tackle the above criticisms, we improved Algorithm 1 by increasing its robustness through several strategies. We propose in Algorithm 2 the enhanced version of Algorithm 1. First, we introduced a relaxation parameter $\alpha \in (0, 1]$ in the fixed-point iteration, by rewriting it as $\mathbf{x}_0^{(k)} \leftarrow \mathbf{x}_0^{(k-1)} + \alpha(\widetilde{\mathbf{x}} - \mathbf{x}^{(k-1)})$. We update $\alpha$ adaptively as written in Algorithm 3: in case of non-convergence of the Newton iterations, we repeat the last fixed-point iteration for a smaller value of $\alpha$; in case of convergence, we increase it in the next iteration. Moreover, the fixed-point iterations are nested inside an outer loop, in which we progressively increase the value of $p_{\mathrm{LV}}$ and $T_{\mathrm{a}}$ until the target value $\widetilde{p}$ and $\widetilde{T}_{\mathrm{a}}$ is reached (continuation method). In this outer loop, we adaptively change the step $\omega^{(k)}$, by decreasing it in case of failure of the inner fixed-point loop and increasing it in case of success.

We remark that after each failure of either SteadyStateMechanics or FixedPoint function, we reset the value of $\mathbf{d}$ to the last solution of Eq. (18). This ensures the success of such continuation strategy, by which we move in the space $(\mathbf{x}_0, \mathbf{d}, \omega)$ staying close to the intersection of the curves $\mathbf{d} = \widetilde{\mathbf{x}} - \mathbf{x}_0$ and $\mathbf{d} = \mathbf{d}_{\mathrm{eq}}(\mathbf{x}_0, \omega p_{\mathrm{LV}}, \omega T_{\mathrm{a}})$ (see Fig. 4).

## 4.2 Projection of the reference configuration from a coarser mesh

The procedure to recover the reference configuration (Algorithm 2) implies to numerically solve the elastostatic problem of Eq. (18) multiple times, until the fixed point algorithm converges. According to our experience, this procedure can be very computational demanding in cardiac applications, especially if realistic human heart geometries and fine representations of the computational domain need to be considered. To overcome this issue, the reference configuration recovery algorithm can be run on a coarser mesh compared to the couple of nested meshes $\mathcal{T}_{\mathrm{h}_2}$-$\mathcal{T}_{\mathrm{h}_1}$ used for the electromechanical model. A natural choice can be to exploit once again the efficient hierarchical octree structure, but starting from a coarser level $\mathrm{h}_3$, generating accordingly the triad of nested meshes $\mathcal{T}_{\mathrm{h}_3}$-$\mathcal{T}_{\mathrm{h}_2}$-$\mathcal{T}_{\mathrm{h}_1}$. However, a further level of coarsening would imply a loss of geometric accuracy that is propagated also into the finer couple of meshes $\mathcal{T}_{\mathrm{h}_2}$-$\mathcal{T}_{\mathrm{h}_1}$, affecting the electromechanical model.

**Algorithm 2** Reference configuration recovery (enhanced version)

```
parameters: k_max, ε_tol^ramp, ε_tol^final, γ_ω^+, γ_ω^−, Δω_max
output: converged_RC, x_0

procedure ReferenceConfiguration(x̃, p̃, T̃_a)
    x_0^(0) ← x̃
    ω^(0) ← 0
    Δω ← Δω_max
    for k = 0, ..., k_max do
        ω^(k) ← min(ω^(k−1) + Δω, 1)
        if ω^(k) = 1 then
            ε_tol ← ε_tol^final
        else
            ε_tol ← ε_tol^ramp
        end if
        (converged_FP, x_0^(k)) ← FixedPoint(x̃, ω^(k) p̃, ω^(k) T̃_a x_0^(k−1), ε_tol)
        if converged_FP then
            if ω^(k) = 1 then
                return (true, x_0^(k))                          ▷ Ramp converged.
            end if
            Δω ← min(γ_ω^+ Δω, Δω_max)
        else
            Δω ← γ_ω^− Δω
        end if
    end for
    return (false, 0)                 ▷ Maximum number of iterations reached.
end procedure
```

**Algorithm 3** Inner fixed-point loop of the reference configuration recovery algorithm

**parameters:** $k_{\max}, \alpha_{\min}, \alpha_{\max}, \gamma_\alpha^+, \gamma_\alpha^-$
**output:** `converged_FP`, $\mathbf{x}_0$

**procedure** FixedPoint($\widetilde{\mathbf{x}}, p_{\mathrm{LV}}, T_{\mathrm{a}}, \mathbf{x}_0, \epsilon_{\mathrm{tol}}$)
    $\mathbf{x}_0^{(0)} \leftarrow \mathbf{x}_0$
    $\alpha \leftarrow \alpha_{\max}$
    (`converged_SSM`, $\mathbf{d}^{(0)}$) $\leftarrow$ SteadyStateMechanics($\mathbf{x}_0^{(0)}, p_{\mathrm{LV}}, T_{\mathrm{a}}$)
    **if not** `converged_SSM` **then**
        **return** (false, $\mathbf{0}$)
    **end if**
    **for** $k = 0, \ldots, k_{\max}$ **do**
        $\mathbf{x}_0^{(k)} \leftarrow \mathbf{x}_0^{(k-1)} + \alpha(\widetilde{\mathbf{x}} - \mathbf{x}^{(k-1)})$ ▷ Fixed-point update.
        (`converged_SSM`, $\mathbf{d}^{(k)}$) $\leftarrow$ SteadyStateMechanics($\mathbf{x}_0^{(k)}, \widetilde{p}, \widetilde{T}_{\mathrm{a}}$)
        **if** `converged_SSM` **then**
            $\mathbf{x}^{(k)} \leftarrow \mathbf{x}_0^{(k)} + \mathbf{d}^{(k)}$
            **if** $\|\mathbf{x}^{(k)} - \widetilde{\mathbf{x}}\| \leq \epsilon_{\mathrm{tol}} \|\mathbf{d}^{(k)}\|$ **then**
                **return** (true, $\mathbf{x}_0^{(k)}$) ▷ Fixed-point converged.
            **end if**
            $\alpha \leftarrow \min(\gamma_\alpha^+ \alpha, \alpha_{\max})$
        **else**
            $\alpha \leftarrow \max(\gamma_\alpha^- \alpha, \alpha_{\min})$
        **end if**
    **end for**
    **return** (false, $\mathbf{0}$) ▷ Maximum number of iterations reached.
**end procedure**

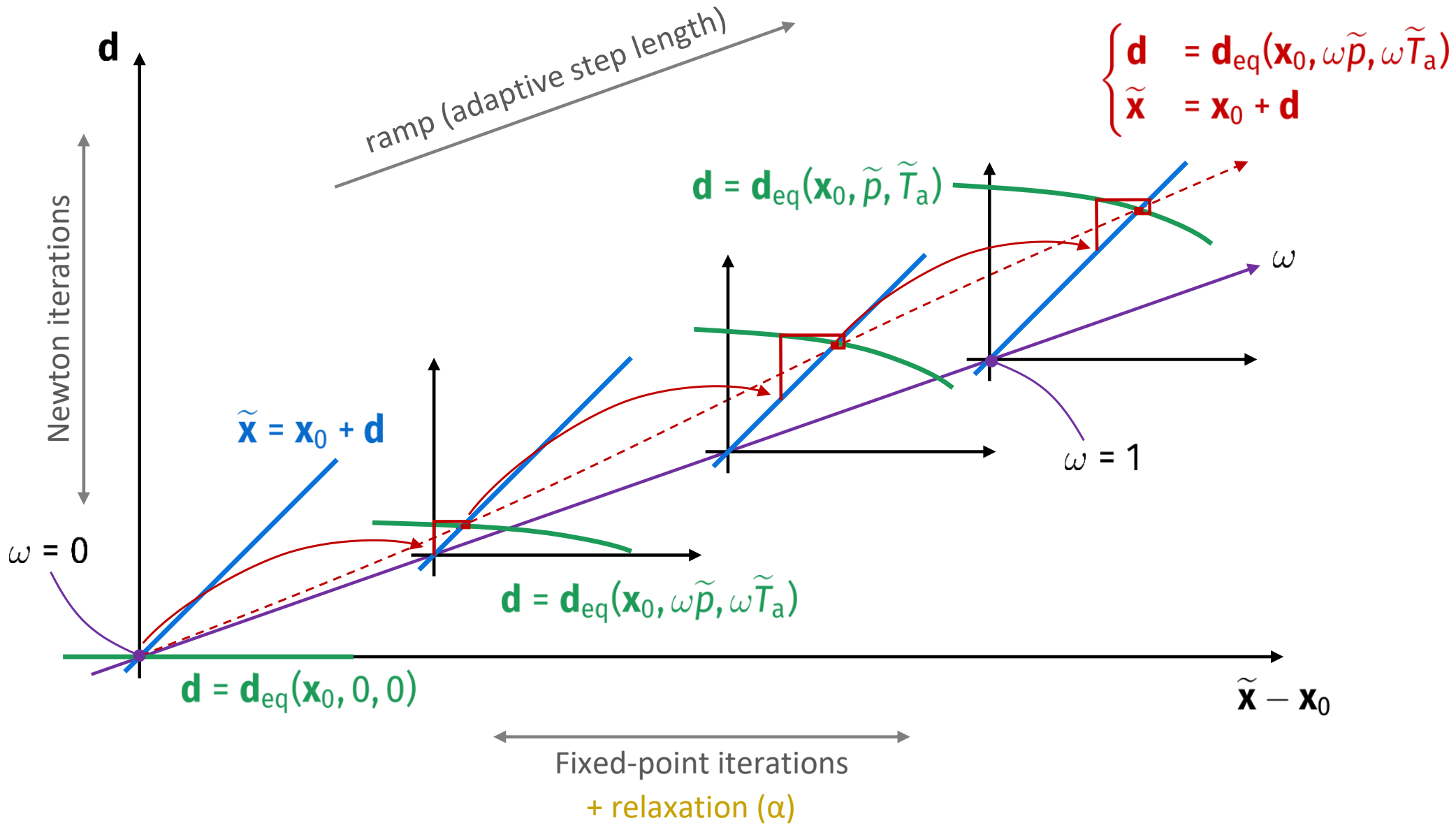

Figure 4: Representation of the enhanced version of the reference configuration recovery algorithm (see Algorithms 2 and 3).

For this reason, we propose a projection technique which enables to map to the reference configuration $\Omega_0$ from a coarser non-nested mesh. This mesh - named $\mathcal{T}_{\widetilde{h}_3}$ - can be independently generated with a mesh size $\widetilde{h}_3$ such that $h_3 < \widetilde{h}_3 < h_2$. Such strategy provides the advantages of having more flexibility on the choice of the mesh size $\widetilde{h}_3$ and of preserving the geometric accuracy for the electromechanical meshes $\mathcal{T}_{h_2}$-$\mathcal{T}_{h_1}$.

The complete procedure consists of the following steps:

1. generate two non-nested computational meshes $\mathcal{T}_{h_2}$ and $\mathcal{T}_{\widetilde{h}_3}$ from the deformed configuration $\widetilde{\Omega}$, i.e. the one reconstructed from the medical images, such that $h_2 < \widetilde{h}_3$. The former is characterized by the target mesh size for the mechanical simulation.
2. solve the reference configuration recovery (Algorithm 2) on the coarser mesh $\mathcal{T}_{\widetilde{h}_3}$, obtaining the displacement field $\mathbf{d}_{\widetilde{h}_3}$;
3. project the displacement $\mathbf{d}_{\widetilde{h}_3}$ on the finer mesh $\mathcal{T}_{h_2}$ obtaining the field $\widehat{\mathbf{d}}_{h_2}$, which approximates the displacement $\mathbf{d}_{h_2}$, i.e. the one that would be computed if the reference configuration recovery (Algorithm 2) is applied directly on the finer mesh $\mathcal{T}_{h_2}$;
4. move each vertex of the mesh $\mathcal{T}_{h_2}$ according to $\mathbf{x}_0 = \widetilde{\mathbf{x}} - \widehat{\mathbf{d}}_{h_2}$, recovering the computational mesh $\mathcal{T}_{h_2}$ that describes the reference configuration $\Omega_0$;
5. hierarchically refine $\mathcal{T}_{h_2}$ to generate the fine mesh $\mathcal{T}_{h_1}$ for the electrophysiology.

We remark that the projection (step 3) is necessary despite both $\mathcal{T}_{h_2}$ and $\mathcal{T}_{\widetilde{h}_3}$ describe the same domain $\widetilde{\Omega}$. Indeed, in practice, their boundaries do not match, since they are independent polygonal surfaces made of piecewise linear elements. Thus, some vertexes of $\mathcal{T}_{h_2}$ can lie outside $\mathcal{T}_{\widetilde{h}_3}$, as illustrated in the 2D sketch of Fig. 5. As a consequence, to recover $\widehat{\mathbf{d}}_{\mathrm{h}_2}$ in all the vertexes of the mesh $\mathcal{T}_{\mathrm{h}_2}$, we proceed in a different way if these vertexes lie inside or outside $\mathcal{T}_{\widetilde{\mathrm{h}}_3}$. In particular:

- for the internal vertexes – denoted as $\widetilde{\mathbf{x}}^i_{\mathrm{h}_2}$, $i = 0, 1, 2, \dots$ – it is sufficient to find the element $K^i \in \mathcal{T}_{\mathrm{h}_3}$ such that $K^i \ni \widetilde{\mathbf{x}}^i_{\mathrm{h}_2}$ and evaluate $\mathbf{d}_{\mathrm{h}_3}(\widetilde{\mathbf{x}}^i_{\mathrm{h}_2})$ exploiting its Finite Element expansion on $K^i$ (see Fig. 5, green element);
- conversely, for each external point – denoted as $\widetilde{\mathbf{x}}^j_{\mathrm{h}_2}$, $j = 0, 1, 2, \dots$ – more sub-steps are necessary, (see Fig. 5, purple element):
  1. we find the closest element $K^j \in \mathcal{T}_{\widetilde{\mathrm{h}}_3}$ from the external point $\widetilde{\mathbf{x}}^j_{\mathrm{h}_2}$;
  2. on $K^j$ we find the closest point $\widetilde{\mathbf{x}}^j_{\mathrm{h}_3}$ to $\widetilde{\mathbf{x}}^j_{\mathrm{h}_2}$;
  3. we evaluate $\mathbf{d}_{\mathrm{h}_3}(\widetilde{\mathbf{x}}^j_{\mathrm{h}_3})$ projecting the resulting value into the external point $\widetilde{\mathbf{x}}^j_{\mathrm{h}_2}$.

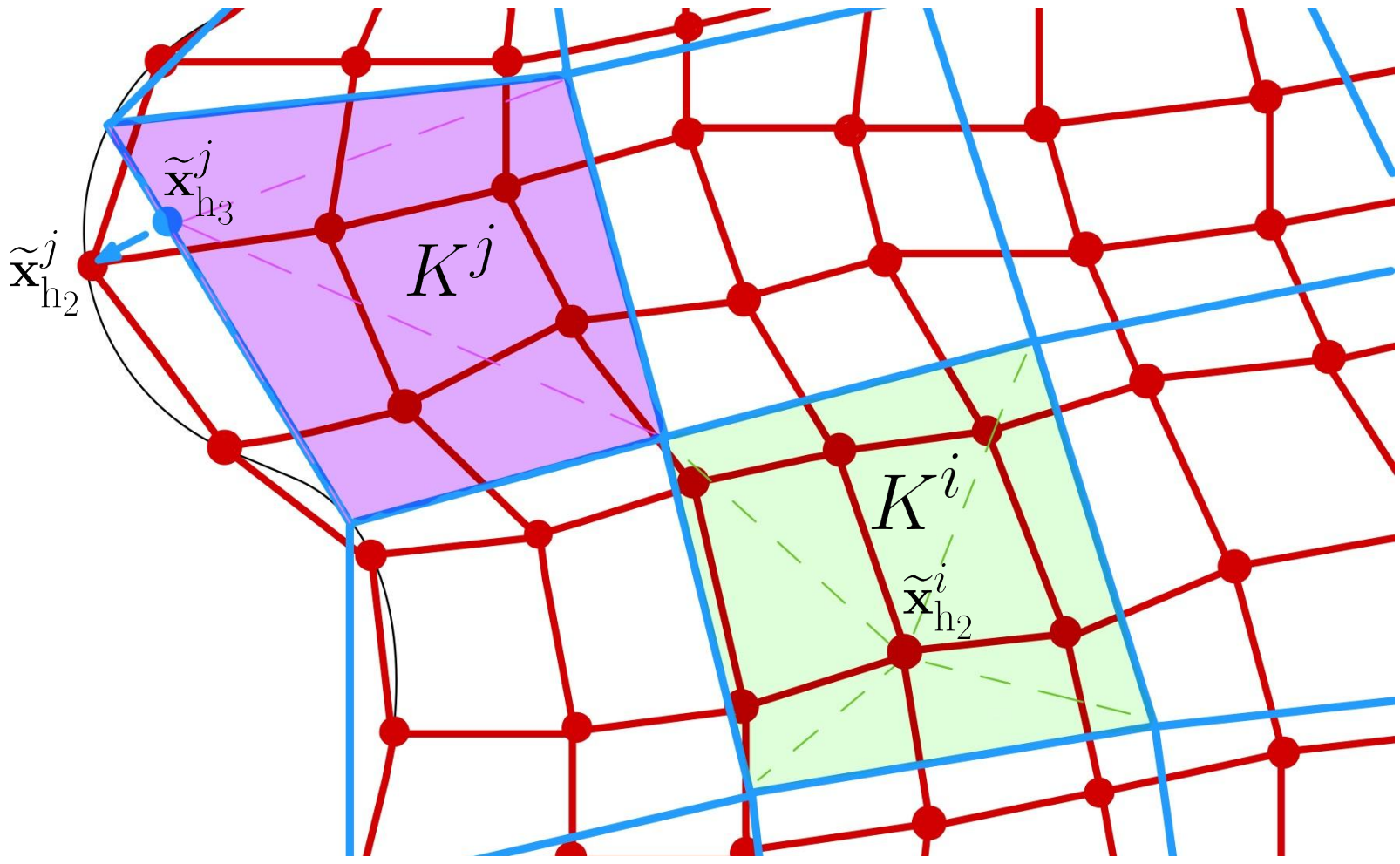

Figure 5: A sketch of the projection procedure from a coarse mesh $\mathcal{T}_{\widetilde{\mathrm{h}}_3}$ (in blue) into a finer mesh $\mathcal{T}_{\mathrm{h}_2}$ (in red): on the internal points $\widetilde{\mathbf{x}}^i_{\mathrm{h}_2} \in \mathcal{T}_{\mathrm{h}_2}$ we recover the value exploiting the basis functions of the element $K^i \in \mathcal{T}_{\mathrm{h}_3}$ (in green); on the external points $\widetilde{\mathbf{x}}^j_{\mathrm{h}_2} \in \mathcal{T}_{\mathrm{h}_2}$ we project the value of the closest point $\widetilde{\mathbf{x}}^j_{\mathrm{h}_3} \in \mathcal{T}_{\widetilde{\mathrm{h}}_3}$ evaluated on the closest element $K^j \in \mathcal{T}_{\widetilde{\mathrm{h}}_3}$ (in purple).

From the implementation point of view, this projection is performed by exploiting the *VTK* library [53] in our Finite Element library `life`$^{\mathtt{x}}$ (`https://lifex.gitlab.io/lifex`). *VTK* filtering utilities allow to locate all internal points, which are the majority, in a really fast way. Moreover, the *VTK* library is efficient in performing closest points interpolation, leading to a very fast projection procedure. However, we remark that the hierarchical octree structure [1, 11] still remains more effective for electromechanical simulations, where the exchange of information between (nested) meshes occurs at each time step. Indeed, the projection presented here is intended to be a single pre-processing step to be performed before the non-stationary electromechanical simulation. An example of this projection in a left ventricle is shown in Sect. 5.1.

| | Number of elements | Number of vertices | $h_{\text{mean}}$ |
|---|---|---|---|
| **Mesh1** | | | |
| Electrophysiology | 1'008'896 | 1'046'641 | 1 mm |
| Activation/Mechanics | 15'764 | 19'099 | 4 mm |
| **Mesh2** | | | |
| Electrophysiology | 2'590'464 | 2'663'817 | 0.75 mm |
| Activation/Mechanics | 40'476 | 47'529 | 3 mm |
| **Mesh3** | | | |
| Electrophysiology | 8'939'776 | 9'116'741 | 0.5 mm |
| Activation/Mechanics | 139'684 | 159'149 | 2 mm |

Table 1: Details of the Zygote left ventricle meshes used for the numerical results (Sec. 5). For each of the three configuration (Mesh1, Mesh2, and Mesh3) two nested meshes are generated: the finer - used for the Electrophysiology model - is obtained by two recursive splitting of each element of the coarser - adopted for both the Activation and Mechanics models.

# 5 Numerical results

We present the numerical simulations of left ventricle electromechanics. The geometry is preprocessed from the Zygote Solid 3D heart model [30], which represents the $50^{\text{th}}$ percentile of a healthy caucasian male in the U.S., reconstructed from an high resolution computed tomography scan. In particular, in Tab. 1 we list the three configurations of the computational meshes under consideration, from the coarsest (Mesh1) to the finest one (Mesh3). For all the cases, we employ for the electrophysiology problem a mesh size four times smaller than for the mechanical problem. The listed meshes were generated using `vmtk` (`www.vmtk.org`) [2], in particular by exploiting the recently proposed tools for cardiac mesh generation [18].

First, we show an example of projection of the reference configuration from a coarser mesh, exploiting the method detailed in Sec. 4.2. Then, we provide a grid refinement study, by considering increasingly refined meshes, employing in all cases different grids for electrophysiology and activation/mechanics. We depict activation maps generated using either $\mathbb{Q}_1$ or $\mathbb{Q}_2$ elements for electrophysiology and only $\mathbb{Q}_1$ element for activation and mechanics: in this way we show that our intergrid transfer operator can handle different Finite Element spaces on different core models, permitting to explore high-order methods in future works. Finally, we show that our model can generate several scenarios according to the chosen parameter sets: specifically, we vary preload, afterload and contractility, and we evaluate the effects of these changes on left ventricle PV loops.

The numerical methods of Sec. 3 have been implemented in `life`$^{\text{x}}$ (`https://lifex.gitlab.io/lifex`), a high-performance `C++` library developed within the iHEART project[1] and based on the `deal.II` (`https://www.dealii.org`) Finite Element core [3]. The numerical simulations have been obtained by employing `life`$^{\text{x}}$ in a parallel setting, using either a HPC resource available at MOX (48 Intel Xeon ES-2640 CPUs) or the GALILEO supercomputer from Cineca (7 nodes endowed with 36 Intel Xeon E5-2697 v4 2.30 GHz CPUs, for a total number of 252 cores).

In terms of the settings our the numerical simulations, we apply a current $\widetilde{\mathcal{I}}_{\text{app}}(\boldsymbol{x}, t)$, distributed in space as a Gaussian with peak $\widetilde{\mathcal{I}}_{\text{app}}^{\max}$, for a duration of $t_{\text{app}}$, in three different regions of the

[1] iHEART - An Integrated Heart Model for the simulation of the cardiac function, European Research Council (ERC) grant agreement No 740132, P.I. Prof. A. Quarteroni

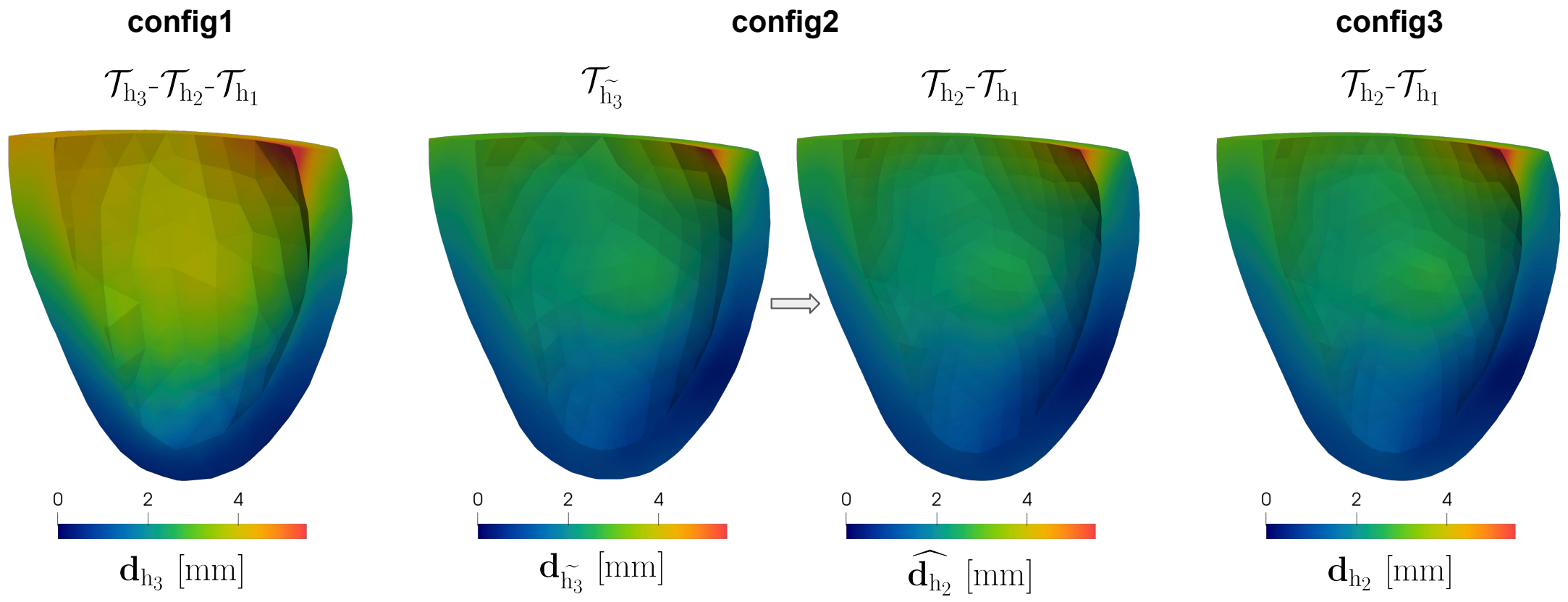

Figure 6: Solving the reference configuration recovery (Algorithm 2) for three different settings, with a common target mesh size for the mechanical simulation ($\mathcal{T}_{\mathrm{h}_2}, h = 3\,\mathrm{mm}$): in *config1* we apply the recovery algorithm on the coarsest mesh of a triad of hierarchical nested meshes ($\mathcal{T}_{\mathrm{h}_3}, h = 6\,\mathrm{mm}$); in *config2* the algorithm runs on an independent coarse mesh ($\mathcal{T}_{\widetilde{\mathrm{h}}_3}, h = 4\,\mathrm{mm}$) and the displacement is projected onto $\mathcal{T}_{\mathrm{h}_2}$; in *config3* we compute the reference configuration directly on $\mathcal{T}_{\mathrm{h}_2}$.

myocardium to trigger the electrical signal in the left ventricle. These regions are located at the central distance between the apex and the base of the ventricle. Even if we do not model the Purkinje network [51, 57], such electrical activation of the tissue that we propose is known to provide meaningful results [21]. We use the Bayer-Blake-Plank-Trayanova algorithm [6, 41] to generate the fibers distribution (field $\mathbf{f}_0$) for our geometry, using $\alpha_{\mathrm{epi}} = -60°$, $\alpha_{\mathrm{endo}} = 60°$, $\beta_{\mathrm{epi}} = 20°$ and $\beta_{\mathrm{endo}} = -20°$. In all numerical tests, we set for simplicity $p_{\mathrm{EX}}(t) = 0$. All the parameters related to the monodomain equation and passive mechanics are reported in Tabs. 2 and 3 respectively. We use a time step $\Delta t_1 = 50\,\mu\mathrm{s}$ for electrophysiology, activation and circulation, whereas mechanics is solved in time with a time step $\Delta t_2 = 250\,\mu\mathrm{s}$. Further information about the parameters of both linear and nonlinear solvers are reported in Appendix B.

## 5.1 Projection of the reference configuration

We discuss the results of the linear projection of the reference configuration (Sect. 4.2) by running the recovery of reference configuration (Algorithm 2) for the Zygote left ventricle geometry. In particular, we consider the following mesh configurations:

- *config1*: a triad of nested meshes $\mathcal{T}_{\mathrm{h}_3}$-$\mathcal{T}_{\mathrm{h}_2}$-$\mathcal{T}_{\mathrm{h}_1}$; the coarsest one for the recovery algorithm, the intermediate one for mechanics and the finer one for electrophysiology;
- *config2*: an independent coarser mesh $\mathcal{T}_{\widetilde{\mathrm{h}}_3}$ for the recovery algorithm, and a couple of nested meshes $\mathcal{T}_{\mathrm{h}_2}$-$\mathcal{T}_{\mathrm{h}_1}$ for the electromechanical model;
- *config3*: a couple of nested meshes $\mathcal{T}_{\mathrm{h}_2}$-$\mathcal{T}_{\mathrm{h}_1}$ for the electromechanical model, where the mesh $\mathcal{T}_{\mathrm{h}_2}$ is used also for the recovery algorithm.

In Fig. 6, we show the results of the recovery algorithm on these three configurations, choosing $h = 3\,\mathrm{mm}$ as target mesh size for mechanics (i.e. $\mathcal{T}_{\mathrm{h}_2}$). More in detail, we create an ad-hoc hierarchical triad $\mathcal{T}_{\mathrm{h}_3}$ ($h = 6\,\mathrm{mm}$)-$\mathcal{T}_{\mathrm{h}_2}$($h = 3\,\mathrm{mm}$)-$\mathcal{T}_{\mathrm{h}_1}$($h = 0.75\,\mathrm{mm}$) for *config1*. Instead, *config2* and *config3* exploit some meshes listed in Tab. 1: Mesh2 is selected as the couple $\mathcal{T}_{\mathrm{h}_2}$-$\mathcal{T}_{\mathrm{h}_1}$, and Mesh1, activation/mechanics as $\mathcal{T}_{\widetilde{\mathrm{h}}_3}$.

The displacement field $\mathbf{d}_{\mathrm{h}_2}$ of *config3* represents the reference result since it is obtained by running directly the recovery algorithm on $\mathcal{T}_{\mathrm{h}_2}$, i.e. the finer mesh that we consider. On the contrary, *config1* presents a large loss of geometric accuracy, which is propagated into the two meshes $\mathcal{T}_{\mathrm{h}_2}$-$\mathcal{T}_{\mathrm{h}_1}$ of the electromechanical model. Indeed, $\mathcal{T}_{\mathrm{h}_3}$ appears too sharp and not realistic, especially at the base and at the apex. Moreover, the displacement field $\mathbf{d}_{\mathrm{h}_3}$ computed on $\mathcal{T}_{\mathrm{h}_3}$ is qualitative different from the reference one ($\mathbf{d}_{\mathrm{h}_2}$, *config3*).

In *config2*, we first compute the displacement field $\mathbf{d}_{\widetilde{\mathrm{h}}_3}$ on the independent coarse mesh $\mathcal{T}_{\widetilde{\mathrm{h}}_3}$, and then we obtain the projected field $\widehat{\mathbf{d}}_{\mathrm{h}_2}$ on $\mathcal{T}_{\mathrm{h}_2}$ by exploiting the linear projection procedure (Sect. 4.2). The projected field $\widehat{\mathbf{d}}_{\mathrm{h}_2}$ features a good qualitative match with the field $\mathbf{d}_{\mathrm{h}_2}$. More quantitatively, we obtain an average error $|\mathbf{d}_{\mathrm{h}_2} - \widehat{\mathbf{d}}_{\mathrm{h}_2}| = 0.105\,\mathrm{mm}$. This quantity is lower than the standard resolution of cardiac medical images (about 0.5 mm for high-quality acquisition). Thus, we can consider this error lower than the geometric uncertainty when dealing with patient-specific simulations. Moreover, being $\mathcal{T}_{\widetilde{\mathrm{h}}_3}$ made of a considerable lower number of elements if compared to $\mathcal{T}_{\mathrm{h}_2}$ (see Tab. 1), *config2* saves a large amount of computational time for the execution of Algorithm 2, compared to *config3*.

*Config2* is the cheapest procedure in terms of total computational cost, moreover it presents a negligible error and does not affect the geometrical accuracy of the electromechanical model. We conclude by recommending this setting when a fine patient-specific mesh is employed for the mechanical problem.

## 5.2 Mesh convergence

To study the results of our numerical model when the space discretization is refined, we run simulations on Mesh1, Mesh2 and Mesh3 of Tab. 1. In Fig. 7 we display the evolution in time of pressure, solid volume, blood volume, and the PV loop, for the three proposed mesh settings. We show that convergence properties in the electromechanical framework are not only influenced by the average diameter $h_{\mathrm{mean}}$ [4], but also by the effective number of DOFs of the mesh. Indeed, we see that Mesh2 with $h_{\mathrm{mean}} \approx 0.75$ mm for electrophysiology and $h_{\mathrm{mean}} \approx 3$ mm for activation and mechanics already provides a good convergence of both mechanical and fluid dynamics relevant properties for the left ventricle. Under the assumption of properly distributed mesh elements, the number of DOFs itself can potentially be an indicator to evaluate the accuracy of the discretization for cardiac simulations, unifying approximations coming from different mesh elements (e.g., tetrahedra, hexaedra and prisms) and different techniques, such as Finite Element Method, Spectral Element Method or Isogeometric Analysis [10, 40, 42].

In Fig. 8 we depict the activation maps for the Mesh3 setting in Tab. 1. We use either $\mathbb{Q}_1$ or $\mathbb{Q}_2$ Finite Elements for electrophysiology, and $\mathbb{Q}_1$ for both activation and mechanics. We show one simulation with $\mathbb{Q}_2$ Finite Elements for electrophysiology to underline that our mathematical discretization can be extended to high order methods, which are known to be more suitable than standard FEM for wave propagation problems [10, 40]. We observe that the activation map resulting from $\mathbb{Q}_2$ elements has slightly more pronounced anisotropy of the isochrones. Indeed, the additional DOFs provided by $\mathbb{Q}_2$ elements allow to better outline the different roles of conductivities $\sigma_{\mathrm{l}}$, $\sigma_{\mathrm{t}}$ and $\sigma_{\mathrm{n}}$, as shown in [41].

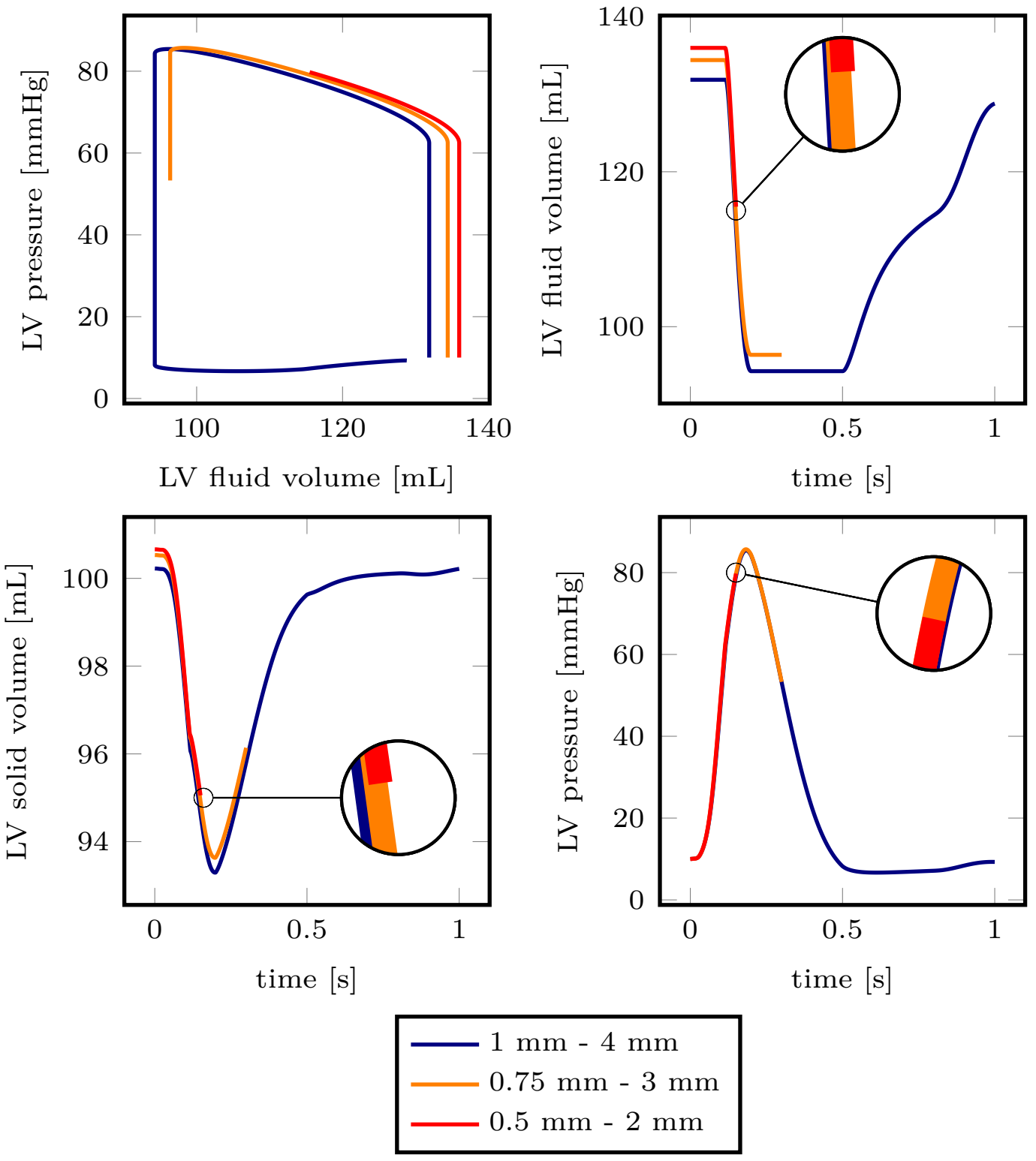

Figure 7: Time evolution of pressure and volumes, and PV loops, related to the Zygote left ventricle, considering different mesh resolutions (Mesh1, Mesh2, Mesh3 of Tab. 1) and $\mathbb{Q}_1$ Finite Element spaces. We report a full heartbeat only for the coarsest mesh (Mesh1).

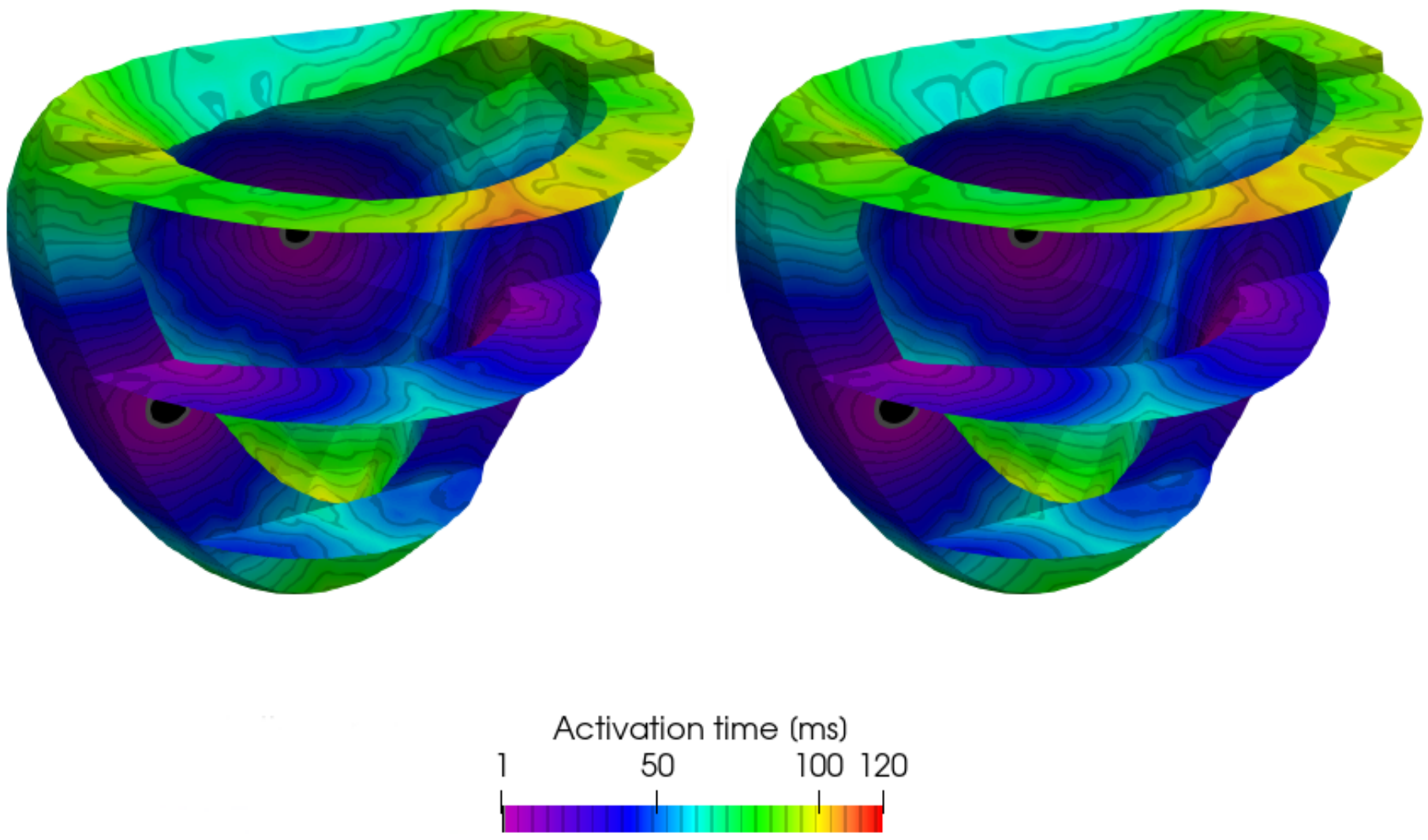

Figure 8: Activation times on the Zygote left ventricle considering different Finite Element spaces ($\mathbb{Q}_1$, 1'046'641 DOFs, on the left, and $\mathbb{Q}_2$, 8'211'745 DOFs, on the right) for electrophysiology.

## 5.3 Tests Cases varying preload, afterload and contractility

We test the response of our electromechanical model to some scenarios of clinical interest. The parameters of the baseline simulation are reported in Appendix A. Starting from this setting, we consider three physiologically relevant scenarios, aimed at investigating the effects of changes in preload, afterload and myocardial contractility, respectively. In all the cases, we simulate several heartbeats, until we reach a periodic regime, but we report only the PV loops related to the last cycle.

### 5.3.1 Baseline simulation

In Fig. 10 we report the time evolution of both the transmembrane potential $u$ and the displacement magnitude $|\mathbf{d}|$ in the Zygote left ventricle, considering the first heartbeat of the simulation. We use mesh configuration 2 in Tab. 1. We compute the initial displacement $\mathbf{d}_0$ by inflating $\Omega_0$ until the desired end-diastolic pressure is reached. In practice we accomplish this task by performing a pressure ramp on the quasi-static approximation of the mechanical problem.

In Fig. 11 we depict the space and time evolutions of the axial stresses $S_{\mathrm{ff}}$, $S_{\mathrm{ss}}$ and $S_{\mathrm{nn}}$, whereas in Fig. 12 we show the shear stresses $S_{\mathrm{fs}}$, $S_{\mathrm{fn}}$ and $S_{\mathrm{sf}}$. We observe that there is no significant generation of stresses (i.e. $S_{\mathrm{ab}} \approx 0$) elsewhere on some specific location at the base, when crossfibers directions are involved. Indeed, the dominant role of fibers and sheets stresses is evident during both the ejection phase (i.e. the second part of systole) and isovolumetric relaxation (i.e. early diastole).

### 5.3.2 Test Case 1

We change the ventricles preload (i.e. their end-diastolic pressure) by modifying the value of the atrial contractility with respect to the baseline setting of Appendix A. More precisely, we consider the two cases when $E_{LA}^{act,max}$ and $E_{RA}^{act,max}$ are respectively increased and decreased by 50%. As shown in Fig. 9, the larger the atrial contractility, the more blood is injected in the ventricle, thus increasing preload. Moreover, our electromechanical model predicts a larger stroke volume for a larger preload. This is consistent with the so-called Frank-Starling effect, a self-regulatory mechanism that guarantees the balance between venous return and cardiac output [31].

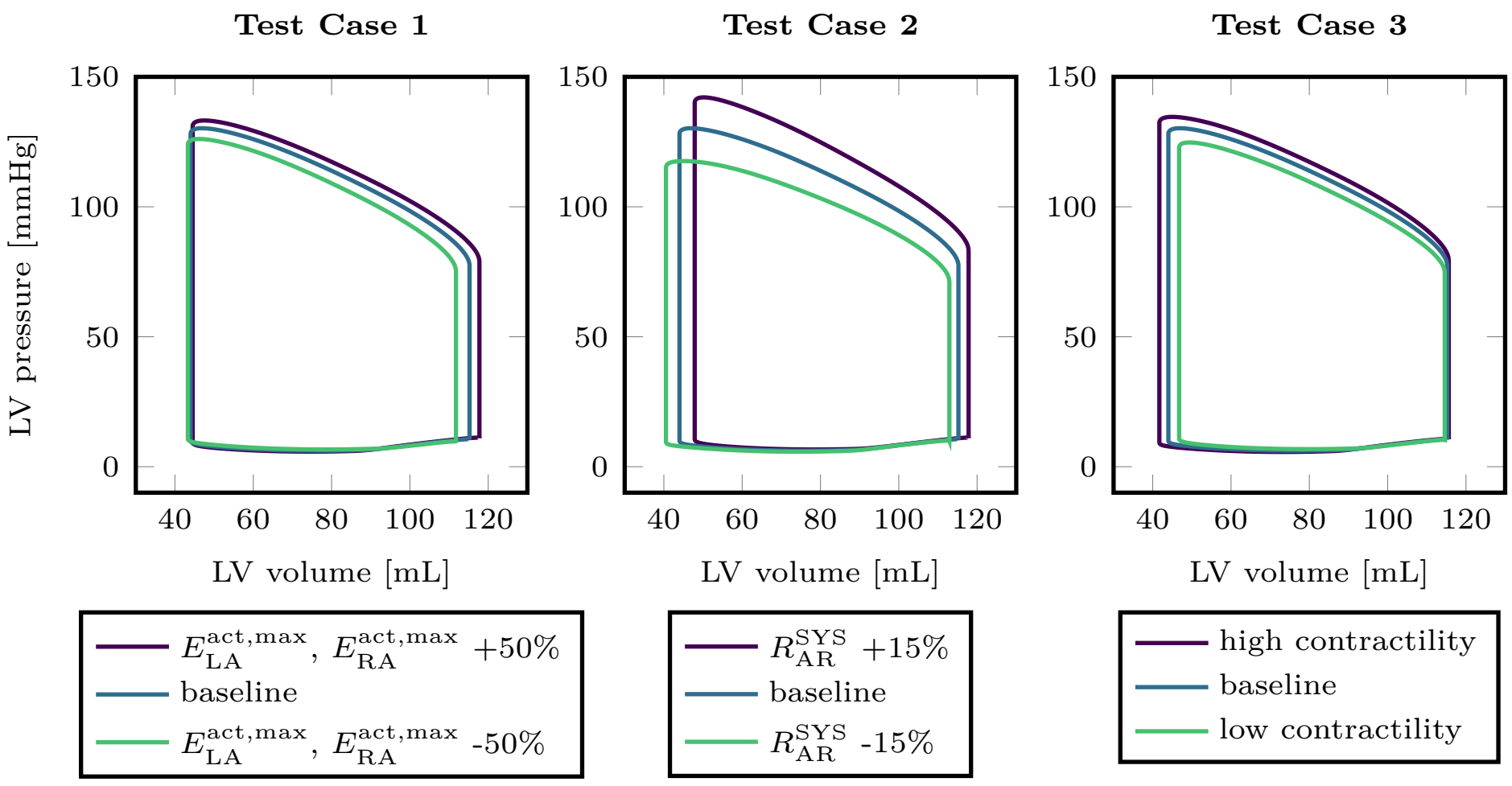

Figure 9: Left ventricle PV loops obtained in the three Test Cases of Sec. 5.3 compared with baseline.

#### 5.3.3 Test Case 2

We investigate the effects of changing the resistance of the arterial circulation. This corresponds to scenarios of a patient affected by diseases associated with hypertension, such as arteriosclerosis, or to the effect of vasodilator-vasoconstrictor drugs. Starting from the baseline setting, we perform two additional simulations where we respectively increase and decrease by 15% the value of $R_{AR}^{SYS}$. In both cases, we modify the value of $C_{AR}^{SYS}$ accordingly, so that the product $R_{AR}^{SYS} C_{AR}^{SYS}$, corresponding to the characteristic time constant of the arterial system, is preserved. The results in Fig. 9 show that an increase of the arterial resistance yields larger values of both the aortic valve opening pressure and the maximal LV pressure (hypertensive effect).

#### 5.3.4 Test Case 3

We consider the response of our electromechanical model to either positive or negative change to inotropic state of the muscle, whose effect is that of increasing and, respectively, decreasing the myocardial contractility. Specifically, starting from the baseline, we first increment and then decrement the atrial contractility by 35% ($E_{LA}^{act,max}$ and $E_{RA}^{act,max}$) by 35% and the ventricular contractility ($T_a^{max}$ and $E_{RV}^{act,max}$). The results in Fig. 9 show that an increase in myocardial contractility generates an increase of both the maximal LV pressure and the stroke volume, but it does not affect preload. As a matter of fact, the end diastolic pressure is unaltered.

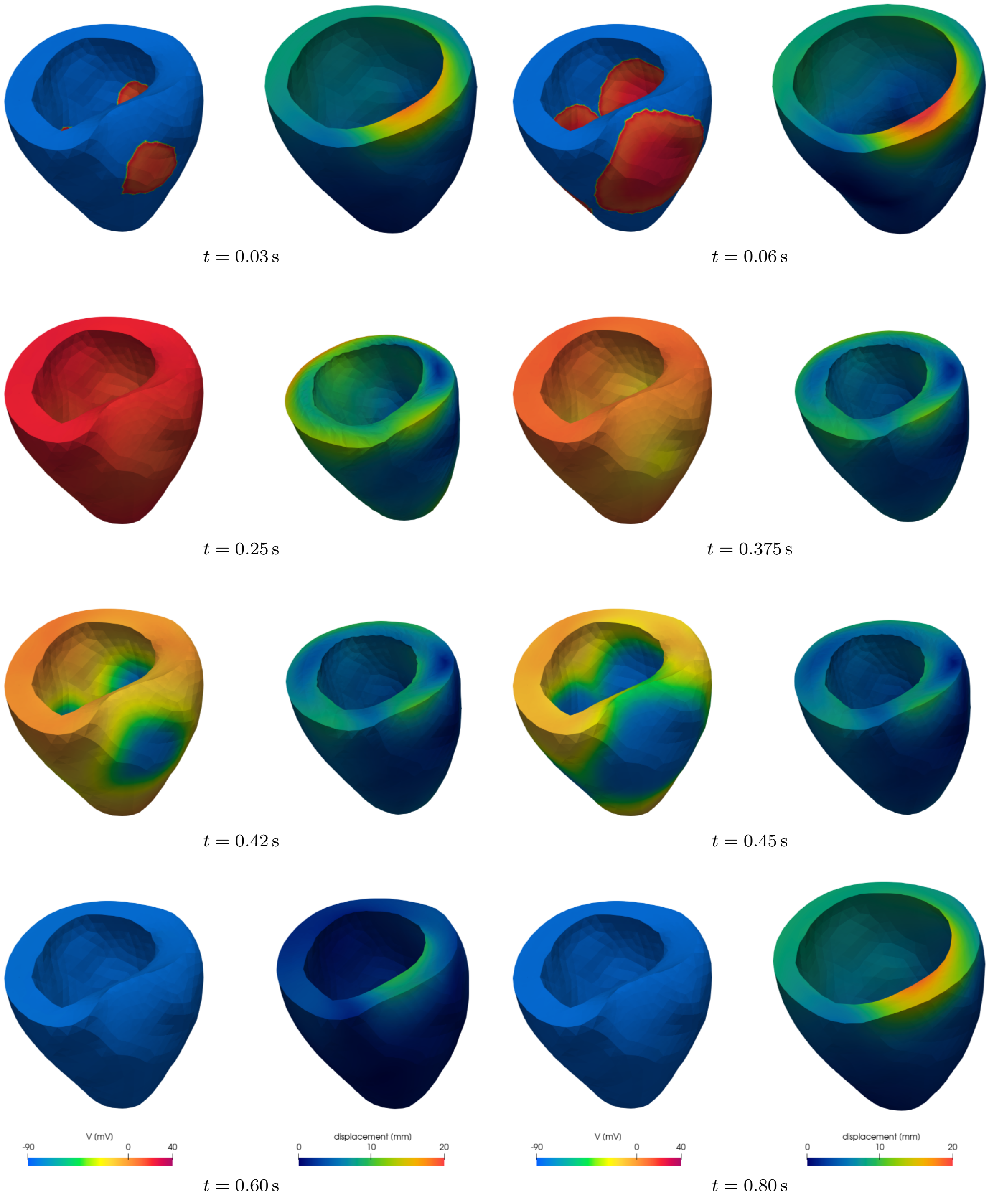

Figure 10: Evolution of the transmembrane potential $u$ and the displacement magnitude $|\mathbf{d}|$ in the Zygote left ventricle over time. The right view of each picture is warped by the displacement vector. Conversely, the potential $u$ is displayed on the reference configuration $\Omega_0$.

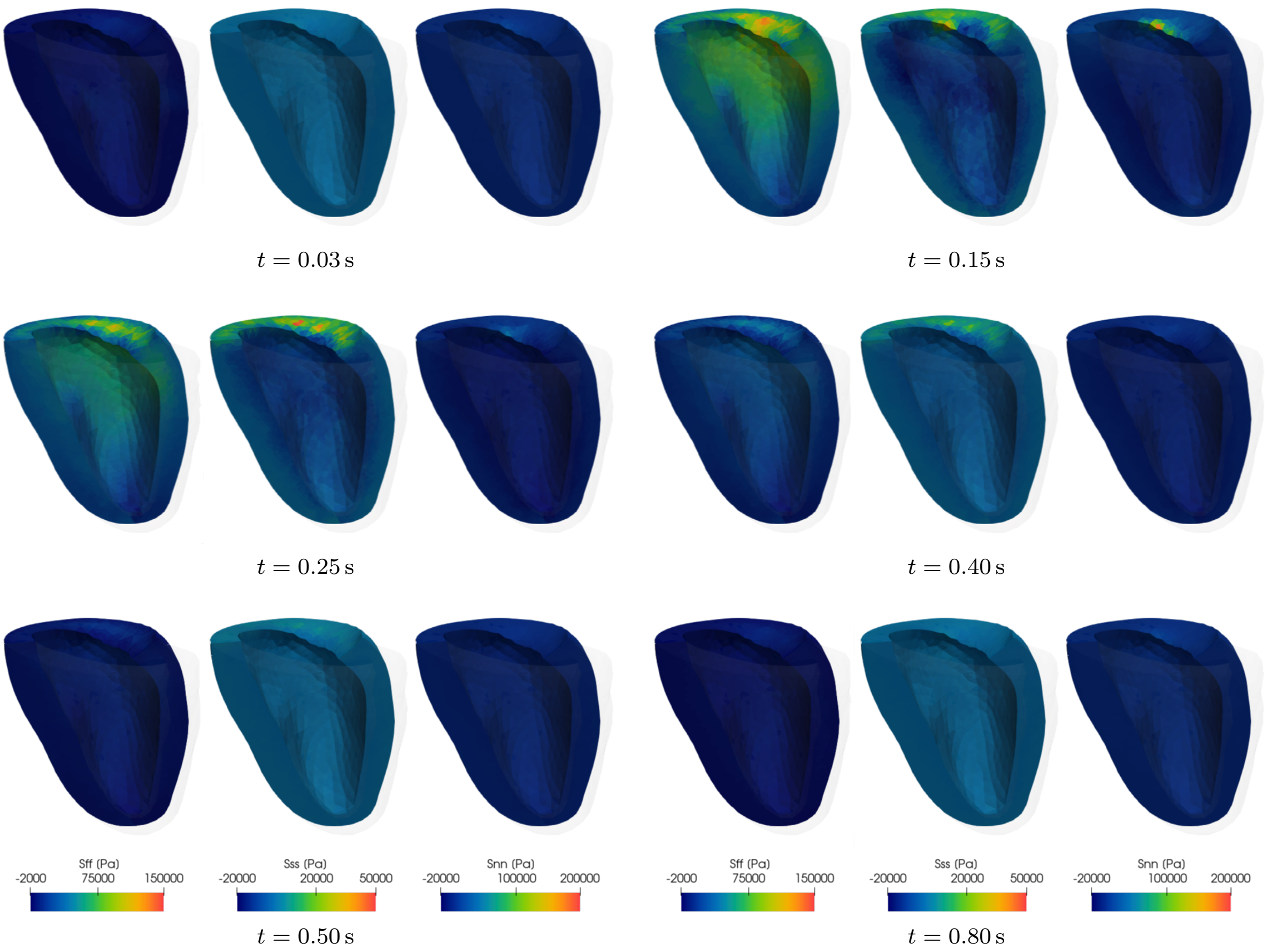

Figure 11: Evolution of $S_{\mathrm{ff}}$, $S_{\mathrm{ss}}$, $S_{\mathrm{nn}}$ during the first heartbeat of the simulation.

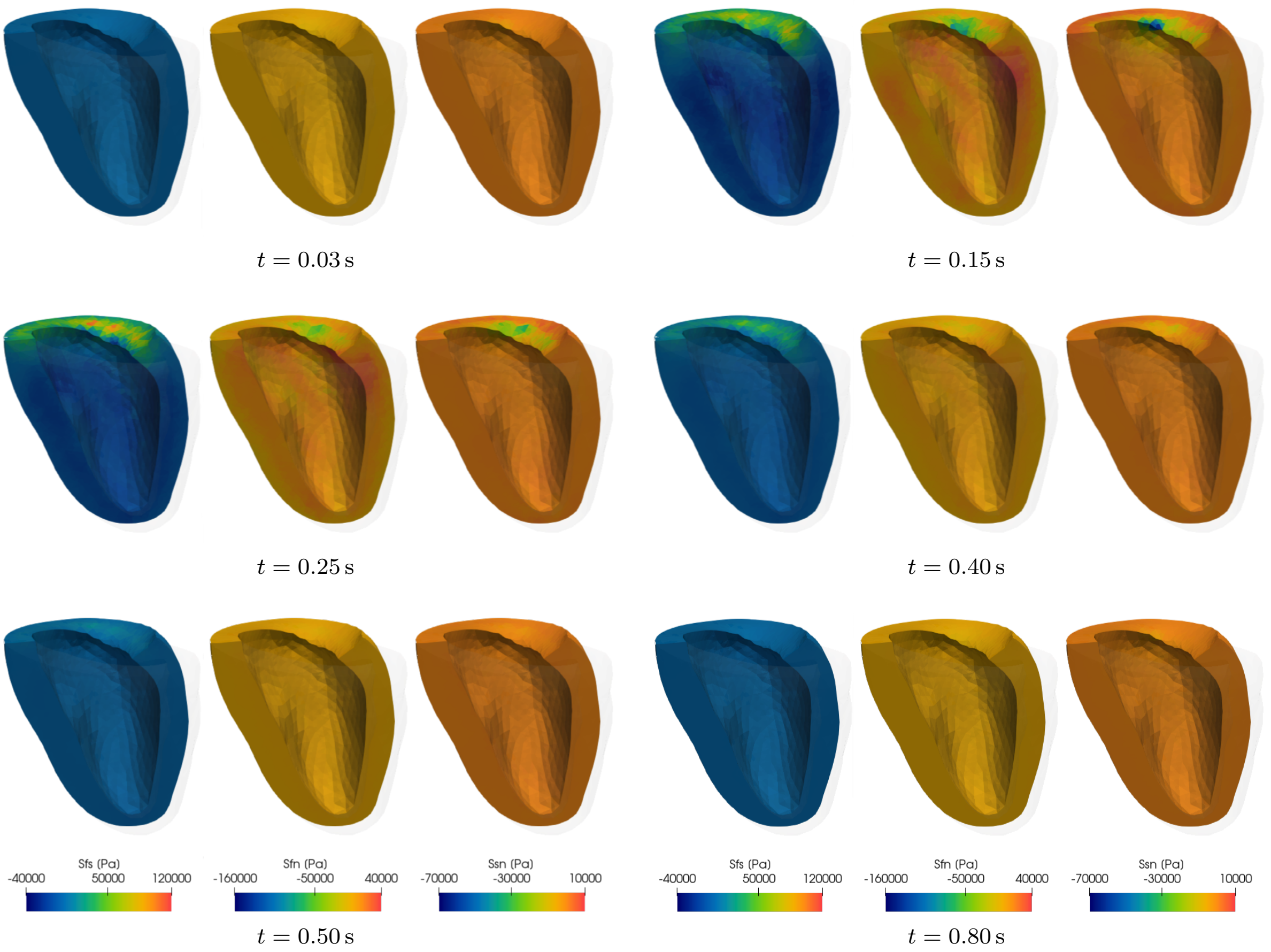

Figure 12: Evolution of $S_{\text{fs}}$, $S_{\text{fn}}$ and $S_{\text{sn}}$ during the first heartbeat of the simulation.

## 6 Conclusions

In this paper, we presented a new numerical model for simulating cardiac electromechanics. The numerical schemes here proposed are aimed at coupling, in a computationally efficient and accurate manner, mathematical models that are significantly biophysically detailed. Our computational framework is based on two main pillars, namely (1) IMEX schemes to approximate the single core models; (2) a fully segregated coupling of the different physics. The advantages yielded by these features are the following.

- We designed IMEX schemes to minimize the number of nonlinear systems that need to be solved, using implicit solvers only for those core models that would otherwise lead to severe CFL restriction on the time step. In particular, our numerical scheme allows to update the variables of the ionic and activation models without the need of solving any algebraic system (neither nonlinear nor linear), and it allows to update the electrical potential by solving a single linear system at each time iteration. The unique nonlinear system of our mathematical model is associated with cardiac mechanics, for which explicit or semi-implicit schemes are unstable, unless a very fine time step is employed, because of to the strong nonlinearities contained in the constitutive law [21].

- We developed a novel scheme to couple the 3D electromechanical model with a 0D circulation model in a fully segregated manner. Our scheme does not present the *balloon dilemma* issue [28] – which affects segregated schemes wherein the displacement update is not aware of the incompressibility constraint of the enclosed fluid, thus possibly leading to the failure of the scheme [5, 25] – by the introduction of a volumetric constraint on the solution of the mechanical problem. In this way, the cavity pressure can be reinterpreted as a Lagrange multiplier associated with the constraint, that enforces the coupling between the 0D circulation model and the 3D electromechanics model. At the algebraic level, we end up with a saddle-point problem, that we solve by means of Schur complement reduction [7].

- The fully segregated scheme proposed in this paper allows to employ different discretizations in space and time to approximate the variables associated with the different physics. This better reflects, at the numerical level, the typical space and time scales characterizing the associated physical phenomena. Specifically, we employ coarser spatial and temporal discretizations for the elastodynamics problem, which is characterized by larger characteristic temporal and spatial scales than the electrophysiological problem [52]. Moreover, it is the most demanding core model from the computational viewpoint, as it involves the solution of a nonlinear system and the time-consuming Jacobian matrix assembly [22].

Regarding the latter aspect, we employ a parallel and flexible intergrid transfer operator [1, 23, 52] that permits to interpolate the solution of a core model among nested meshes and between possibly different FEM spaces [1]. As a matter of fact, among the numerical results presented in this paper, we considered the case where we employ different polynomial orders for the electrophysiology variables, on one hand, and for the activation and mechanics variables, on the other hand. This can be considered as a first step towards possible investigation of high-order methods in our cardiac electromechanics computational framework. Indeed, wave propagation problems, such as the one arising in cardiac electrophysiology, can be more suitably represented using high-order basis functions than classical FEM [10, 40].

As a further novel contribution, we presented an algorithm to reconstruct the reference (i.e. stress-free) configuration of the left ventricle starting from a stressed configuration coming from imaging, by solving a suitable inverse problem. Determining such configuration is essential to correctly initialize electromechanical simulations, which is especially useful in patient-specific scenarios where the end diastolic pressure and/or the end diastolic volume are possibly known.

Finally, we presented several simulations varying different parameters of the model, to affect preload, afterload and contractility, thus investigating the response of our model to scenarios of clinical interest. Our model correctly reproduces the increase of stroke volume as a consequence of increased preload, coherently with the Frank-Starling law [31], thus guaranteeing the matching between the venous return and the cardiac output.

Albeit in this paper we focused on left ventricle simulations, the 3D-0D coupling that we describe (along with its numerical treatment) can be extended to a full four chambers representation of the human heart, a model that we will pursue in forthcoming publications. Similarly, the pipeline to reconstruct the reference geometry can be generalized and employed for both atria and ventricles.

## Acknowledgements

This project has received funding from the European Research Council (ERC) under the European Union's Horizon 2020 research and innovation programme (grant agreement No 740132, iHEART - An Integrated Heart Model for the simulation of the cardiac function, P.I. Prof. A. Quarteroni). We acknowledge the CINECA award under the class C ISCRA project HP10CWQ2GS, for the availability of high performance computing resources and support.

# Appendices

## A Model parameters

We provide the full list of parameters adopted for the simulations referred as baseline in Sec. 5.3. Specifically, Tab. 2 contains the parameters related to the electrophysiological model, Tab. 3 those related to the mechanical model and, finally, Tab. 4 contains the parameters of the circulation model. For the TTP06 model, we adopt the parameters reported in the original paper (for M cells) [55]. For the RDQ18 model, we employ the parameters of the original paper [46].

| Variable | Value | Unit | Variable | Value | Unit |
|---|---|---|---|---|---|
| **Conductivity tensor** | | | **Applied current** | | |
| $\sigma_\mathrm{l}$ | $0.7643 \cdot 10^{-3}$ | $\mathrm{m^2\,s^{-1}}$ | $\widetilde{\mathcal{I}}_{\mathrm{app}}^{\max}$ | 35 | $\mathrm{V\,s^{-1}}$ |
| $\sigma_\mathrm{t}$ | $0.3494 \cdot 10^{-3}$ | $\mathrm{m^2\,s^{-1}}$ | $t_{\mathrm{app}}$ | $3 \cdot 10^{-3}$ | s |
| $\sigma_\mathrm{n}$ | $0.1125 \cdot 10^{-3}$ | $\mathrm{m^2\,s^{-1}}$ | | | |

Table 2: Parameters of the electrophysiological model.

| Variable | Value | Unit | Variable | Value | Unit |
|---|---|---|---|---|---|
| **Constitutive law** | | | **Boundary conditions** | | |
| $B$ | $50 \cdot 10^3$ | Pa | $K_\perp^{\mathrm{epi}}$ | $2 \cdot 10^5$ | $\mathrm{Pa\,m^{-1}}$ |
| $C$ | $0.88 \cdot 10^3$ | Pa | $K_\parallel^{\mathrm{epi}}$ | $2 \cdot 10^4$ | $\mathrm{Pa\,m^{-1}}$ |
| $b_{\mathrm{ff}}$ | 8 | – | $C_\perp^{\mathrm{epi}}$ | $2 \cdot 10^4$ | $\mathrm{Pa\,s\,m^{-1}}$ |
| $b_{\mathrm{ss}}$ | 6 | – | $C_\parallel^{\mathrm{epi}}$ | $2 \cdot 10^3$ | $\mathrm{Pa\,s\,m^{-1}}$ |
| $b_{\mathrm{nn}}$ | 3 | – | | | |
| $b_{\mathrm{fs}}$ | 12 | – | **Activation** | | |
| $b_{\mathrm{fn}}$ | 3 | – | $T_{\mathrm{a}}^{\max}$ | $180 \cdot 10^3$ | Pa |
| $b_{\mathrm{sn}}$ | 3 | – | $SL_0$ | 2 | µm |
| $\rho_{\mathrm{s}}$ | $10^3$ | $\mathrm{kg\,m^{-3}}$ | | | |

Table 3: Parameters of the mechanical model.

| Variable | Value | Unit | Variable | Value | Unit |
|---|---|---|---|---|---|
| **External circulation** | | | **Cardiac chambers** | | |
| $R_{\mathrm{AR}}^{\mathrm{SYS}}$ | 0.8 | $\mathrm{mmHg\,s\,mL^{-1}}$ | $E_{\mathrm{LA}}^{\mathrm{pass}}$ | 0.09 | $\mathrm{mmHg\,mL^{-1}}$ |
| $R_{\mathrm{AR}}^{\mathrm{PUL}}$ | 0.1625 | $\mathrm{mmHg\,s\,mL^{-1}}$ | $E_{\mathrm{RA}}^{\mathrm{pass}}$ | 0.07 | $\mathrm{mmHg\,mL^{-1}}$ |
| $R_{\mathrm{VEN}}^{\mathrm{SYS}}$ | 0.26 | $\mathrm{mmHg\,s\,mL^{-1}}$ | $E_{\mathrm{RV}}^{\mathrm{pass}}$ | 0.05 | $\mathrm{mmHg\,mL^{-1}}$ |
| $R_{\mathrm{VEN}}^{\mathrm{PUL}}$ | 0.1625 | $\mathrm{mmHg\,s\,mL^{-1}}$ | $E_{\mathrm{LA}}^{\mathrm{act,max}}$ | 0.07 | $\mathrm{mmHg\,mL^{-1}}$ |
| $C_{\mathrm{AR}}^{\mathrm{SYS}}$ | 1.2 | $\mathrm{mL\,mmHg^{-1}}$ | $E_{\mathrm{RA}}^{\mathrm{act,max}}$ | 0.06 | $\mathrm{mmHg\,mL^{-1}}$ |
| $C_{\mathrm{AR}}^{\mathrm{PUL}}$ | 10.0 | $\mathrm{mL\,mmHg^{-1}}$ | $E_{\mathrm{RV}}^{\mathrm{act,max}}$ | 0.55 | $\mathrm{mmHg\,mL^{-1}}$ |
| $C_{\mathrm{VEN}}^{\mathrm{SYS}}$ | 60.0 | $\mathrm{mL\,mmHg^{-1}}$ | $V_{0,\mathrm{LA}}$ | 4.0 | mL |
| $C_{\mathrm{VEN}}^{\mathrm{PUL}}$ | 16.0 | $\mathrm{mL\,mmHg^{-1}}$ | $V_{0,\mathrm{RA}}$ | 4.0 | mL |
| $L_{\mathrm{AR}}^{\mathrm{SYS}}$ | $5 \cdot 10^{-3}$ | $\mathrm{mmHg\,s^2\,mL^{-1}}$ | $V_{0,\mathrm{RV}}$ | 10.0 | mL |
| $L_{\mathrm{AR}}^{\mathrm{PUL}}$ | $5 \cdot 10^{-4}$ | $\mathrm{mmHg\,s^2\,mL^{-1}}$ | **Cardiac valves** | | |
| $L_{\mathrm{VEN}}^{\mathrm{SYS}}$ | $5 \cdot 10^{-4}$ | $\mathrm{mmHg\,s^2\,mL^{-1}}$ | $R_{\min}$ | 0.0075 | $\mathrm{mmHg\,s\,mL^{-1}}$ |
| $L_{\mathrm{VEN}}^{\mathrm{PUL}}$ | $5 \cdot 10^{-4}$ | $\mathrm{mmHg\,s^2\,mL^{-1}}$ | $R_{\max}$ | 75006.2 | $\mathrm{mmHg\,s\,mL^{-1}}$ |

Table 4: Parameters of the circulation model. We always consider a heartbeat period $T = 0.8\,\mathrm{s}$.

# B Parameters of the linear and nonlinear solvers

We report the setting used for the linear (Tab. 5) and nonlinear (Tab. 6) solvers to produce the results shown in this paper.

| Physics/Fields | Linear solver | Abs. tol. |
|---|---|---|
| Monodomain model | CG | $10^{-10}$ |
| Activation | GMRES | $10^{-10}$ |
| Mechanics | GMRES | $10^{-8}$ |

Table 5: Tolerances about the linear solver for the different physics.

| Physics/Fields | nonlinear solver | Rel. tol. | Abs. tol. |
|---|---|---|---|
| Mechanics | Quasi-Newton | $10^{-10}$ | $10^{-8}$ |
| Reference configuration | Newton | $10^{-10}$ | $10^{-8}$ |

Table 6: Tolerances about the nonlinear solver for the mechanical problem.